\pdfoutput=1
\documentclass[a4paper, 11pt]{amsart}
\usepackage[utf8]{inputenc}
\usepackage[T1]{fontenc}
\usepackage[american]{babel}
\usepackage{amssymb}
\usepackage{amsmath}
\usepackage{amsthm}
\usepackage{amscd}
\usepackage{amsfonts}
\usepackage{stmaryrd}
\usepackage{pb-diagram}
\usepackage{epic,eepic,epsfig}
\usepackage{a4wide}
\usepackage{nextpage}
\usepackage{fancyhdr}
\usepackage{tikz}
\usepackage{verbatim}
\usepackage{comment}
\usepackage{csquotes}
\usepackage{textcomp}

\usepackage[backref=true,giveninits=true,doi=false,url=false,isbn=false,backend=biber,sorting=nyt,sortcites=true]{biblatex}
    \newbibmacro{string+doiurlisbn}[1]{%
      \iffieldundef{doi}{%
        \iffieldundef{url}{%
          \iffieldundef{isbn}{%
            \iffieldundef{issn}{%
              #1%
            }{%
              \href{http://books.google.com/books?vid=ISSN\thefield{issn}}{#1}%
            }%
          }{%
            \href{http://books.google.com/books?vid=ISBN\thefield{isbn}}{#1}%
          }%
        }{%
          \href{\thefield{url}}{#1}%
        }%
      }{%
        \href{http://dx.doi.org/\thefield{doi}}{#1}%
      }%
    }
    \DeclareFieldFormat{title}{\usebibmacro{string+doiurlisbn}{\mkbibemph{#1}}}
    \DeclareFieldFormat{year}{\usebibmacro{year-parenthesis}{\mkbibemph{#1}}}

\DeclareFieldFormat[article,incollection,thesis,misc,inproceedings,online,inbook]{title}{\usebibmacro{string+doiurlisbn}{\mkbibemph{#1}}}
\AtEveryBibitem{%
  \clearfield{volume}%
  \clearfield{number}}

\usepackage{amsmath}
\usepackage{amsthm}
\usepackage{amssymb} 

\usepackage{amsfonts}

\usepackage{tikz-cd,tikz}
\usetikzlibrary{arrows,backgrounds}

\usepackage{mathtools}
\usepackage{stmaryrd}

\usepackage{hyperref}
\hypersetup{colorlinks=true,allcolors=blue}
\usepackage[capitalize,nameinlink,noabbrev]{cleveref}

\usepackage{enumitem}
\usepackage{wrapfig}

\theoremstyle{plain}
\newtheorem{theorem}{Theorem}[section]

\newtheorem{corollary}[theorem]{Corollary}
\newtheorem{proposition}[theorem]{Proposition}
\newtheorem{conjecture}[theorem]{Conjecture}

\theoremstyle{definition}
\newtheorem{definition}[theorem]{Definition}

\theoremstyle{remark}

\newtheorem{remark}[theorem]{Remark}
\newtheorem{example}[theorem]{Example}

    \newcommand\restr[2]{{
    		\left.\kern-\nulldelimiterspace #1 \right|_{#2}
    }}

    \DeclareFontFamily{U}{wncy}{}
    \DeclareFontShape{U}{wncy}{m}{n}{<->wncyr10}{}
    \DeclareSymbolFont{mcy}{U}{wncy}{m}{n}
    \DeclareMathSymbol{\B}{\mathord}{mcy}{"42} 
    \DeclareMathSymbol{\Sha}{\mathord}{mcy}{"58}

\makeatletter
\newcommand{\subalign}[1]{%
  \vcenter{%
    \Let@ \restore@math@cr \default@tag
    \baselineskip\fontdimen10 \scriptfont\tw@
    \advance\baselineskip\fontdimen12 \scriptfont\tw@
    \lineskip\thr@@\fontdimen8 \scriptfont\thr@@
    \lineskiplimit\lineskip
    \ialign{\hfil$\m@th\scriptstyle##$&$\m@th\scriptstyle{}##$\hfil\crcr
      #1\crcr
    }%
  }%
}
\makeatother

\begin{document}

\title{On the Northcott property for special values of L-functions}
\author{{F}abien {P}azuki}
\author{{R}iccardo {P}engo}
\address{Fabien Pazuki - University of Copenhagen, Department of Mathematical Sciences, Universitetsparken 5, 2100 Copenhagen, Denmark}
\email{\href{mailto:fpazuki@math.ku.dk}{fpazuki@math.ku.dk}}
\address{Riccardo Pengo - Institut für Algebra, Zahlentheorie und Diskrete Mathematik, Fakultät für Mathematik und Physik, Leibniz Universität Hannover, Welfengarten 1, 30167 Hannover, Germany}
    \email{\href{mailto:pengo@math.uni-hannover.de}{pengo@math.uni-hanonver.de}}

\begin{abstract}
We propose an investigation on the Northcott, Bogomolov and Lehmer properties for special values of $L$-functions. We first introduce an axiomatic approach to these three properties. We then focus on the Northcott property for special values of $L$-functions. In the case of $L$-functions of pure motives, we prove a Northcott property for special values located at the left of the critical strip, assuming that the $L$-functions in question satisfy some expected properties. Inside the critical strip, focusing on the Dedekind zeta function of number fields, we prove that such a property does not hold for the special value at one, but holds for the special value at zero, and we give a related quantitative estimate in this case. 
\end{abstract}
\maketitle

\begin{center}
\today
\end{center}

{\flushleft
}

\begin{center}
---------
\end{center}

\thispagestyle{empty}

\section{Introduction}

\subsection{Diophantine properties of the Weil height}
Traditionally, the expression \textit{Northcott property} is used to talk about the finiteness of the set of algebraic numbers having simultaneously bounded height and degree, proved by Northcott (see \cite{no49} and \cite[Theorem~1.6.8]{bg06}). More generally, one can say that a field $F \subseteq \overline{\mathbb{Q}}$ has the \textit{Northcott property} if the sets of elements of $F$ having bounded height are finite. Hence, Northcott's theorem can be reformulated by saying that number fields have the Northcott property. They are not the only fields sharing this property, as we recall in \cref{ex:Weil_height}.

The Northcott property can be relaxed by asking for which fields $F \subseteq \overline{\mathbb{Q}}$ the image of the height $h(F) \subseteq \mathbb{R}_{\geq 0}$ does not admit zero as an accumulation point. If this happens, one says that $F$ has the \textit{Bogomolov property} (introduced in \cite{bz01}), and in recent years there has been an increasing interest in finding fields having the Bogomolov property.
Finally, the Bogomolov property can also be relaxed by looking at fields $F \subseteq \overline{\mathbb{Q}}$ such that the product of the height and the degree of those algebraic numbers which belong to $F$ does not accumulate towards zero. It seems reasonable to us to define this as the \textit{Lehmer property} in view of the famous conjecture of Lehmer which says that $\overline{\mathbb{Q}}$ (and hence each of its subfields) satisfies this property. 
We note that there are many explicit classes of sub-fields of $\overline{\mathbb{Q}}$ which are known to have the Bogomolov property but not the Northcott property. 
It is more difficult to construct fields which provably have the Lehmer property and do not have the Bogomolov one: the first example of a field of this kind, to the best of the authors' knowledge, has been constructed by Amoroso in \cite[Theorem~3.3]{Amoroso_2016}.
We refer the interested reader to the work of the first author, Technau and Widmer (see \cite[Theorem~4]{Pazuki_Technau_Widmer_2022}) for other examples of fields with the Lehmer but not the Bogomolov property, and to \cref{ex:Weil_height} for a more detailed discussion.

\subsection{Heights of motives and special values of \texorpdfstring{$L$}{L}-functions}
We are interested in the Northcott, Bogomolov and Lehmer properties for more general notions of heights, related to objects in algebraic geometry. We will define these properties in \cref{sec:properties_of_heights} for any set-theoretic function $h \colon S \to \Gamma$, where $\Gamma$ is a partially ordered set. In particular, this allows to define such properties for sets of heights, by considering their set-theoretic product (see \cref{def:product_height}). Moreover, this language allows to put within the same framework some results in different areas of Diophantine geometry, as we observe in \cref{sec:examples_minima} and \cref{sec:examples_heights}. 

One can view heights as a way of measuring the complexity of arithmetic or geometric objects, such as algebraic numbers, algebraic number fields, abelian varieties or Galois representations.
More generally speaking, one considers \textit{mixed motives}, which can be thought of as pieces cut out from the cohomology of algebraic varieties defined over a number field $F$.  
A recent pioneering work by Kato \cite{ka18} puts forward several tentative definitions of heights of mixed motives, whose Northcott properties would have momentous consequences, as we recall in \cref{ex:Heights_for_motives}. One way to derive a Northcott property for a height is to compare it to another height, where the statement would be easier to prove. This is a successful strategy for the Faltings height of an abelian variety, which can for instance be explicitly compared to its theta height, as shown in \cite{Pazuki_2012}.
Therefore, it seems natural to try to relate Kato's heights with other kinds of heights for which proving a Northcott property may be easier. 

One of the aims of the present paper is to propose the idea that special values of $L$-functions may serve as such a height. More precisely, if one fixes a prime number $\ell \in \mathbb{N}$, one can associate to every mixed motive $X$ defined over a number field $F$ an $\ell$-adic realisation $R_\ell(X)$, which is a representation of the absolute Galois group $G_F := \operatorname{Gal}(\overline{F}/F)$. Assembling together all the actions of the Frobenius elements belonging to $G_F$, one gets a formal Euler product $L(R_\ell(X),s)$, as explained for instance in \cite{de94}. Then, one expects that this Euler product actually defines a meromorphic function $L(R_\ell(X),s) \colon \mathbb{C} \dashrightarrow \mathbb{C}$, to which one can associate the special values:
\[
    L^\ast(X,n) := \lim_{s \to n} \frac{L(X,s)_\sigma}{(s - n)^{\operatorname{ord}_{s = n}(L(X,s)_\sigma)}}
\]
taken at each integer $n \in \mathbb{Z}$.
Such special values have historically been related to many different kinds of heights. For instance:
\begin{enumerate}
    \item the conjectures of Boyd \cite{bo98} relate certain special values of $L$-functions to the Mahler measure of integral polynomials (see \cref{ex:Mahler_measure});
    \item the formula of Gross and Zagier \cite{Gross_Zagier_1986} relates special values of $L$-functions of modular forms to heights of Heegner points (see \cref{ex:canonical_height});
    \item the conjecture of Colmez \cite{co93} relates logarithmic derivatives of Artin $L$-functions to Faltings's heights (see \cref{ex:faltings_height});
    \item special values of Dedekind $\zeta$-functions appear to be related to the volume of hyperbolic manifolds, which can be regarded as a height (see \cref{ex:volume_hyperbolic_manifolds});
    \item the conjectures of Bloch and Kato can be stated as a relation between special values of $L$-functions and height pairings of algebraic cycles (see \cite{Bloch_Kato_1990}, as well as Fontaine's survey \cite{fo92}).
\end{enumerate}
Thus it seems natural to us to investigate which properties that are typical of heights hold also for special values of $L$-functions. The main results of the present paper go in this direction, as we explain in the next sub-section. Moreover, as we note in \cref{ex:Heights_for_motives}, some ingredients appearing in the definition of Kato's height of motives prompt us to believe that there might be connections between this height and special values of $L$-functions. This will be the subject of future research.

\subsection{Main results}

Let us now summarise the main results of this paper, which concern the Northcott property for special values of $L$-functions. 

First of all, we prove a general result about $L$-functions associated to pure motives, and their special values at the left of the critical strip. More precisely, for every number field $F$ we let $\mathcal{M}_F$ denote the category of pure motives for absolute Hodge cycles, introduced by Deligne in \cite{Deligne_1979}, and for every $w \in \mathbb{Z}$ we denote by $\mathcal{M}^{(w)}(F)$ the set of isomorphism classes of pure objects $X \in \mathcal{M}_F$ of weight $w$, whose $L$-functions satisfy some expected properties, outlined in \cref{def:Mw}. In particular, we suppose that these $L$-functions have the expected region of absolute convergence, and can be meromorphically continued to functions which satisfy the expected functional equation, so that each integer $n \in \mathbb{Z}$ will lie on the left of the critical strip if and only if $2 n < w$.
With this setting in mind, we show that the special values at the left of the critical strip of the $L$-functions associated to pure motives of fixed weight and dimension satisfy the Northcott property, as summarised in the following theorem.

\begin{theorem}[Northcott property at the left of the critical strip]
\label{thm:Northcott_left_critical_strip}
    Let $F$ be a number field.
    Fix an integer $w \in \mathbb{Z}$, and let $\mathcal{M}^{(w)}(F)$ be the set introduced in \cref{def:Mw}.
    Then, for every $\mathbf{B} = (B_1, B_2) \in \mathbb{R}_{\geq 0}^2$ and every $n \in \mathbb{Z}$ such that $2 n < w$, the set
    \begin{equation}\label{eq:finiteMM}
        \mathcal{M}^{(w)}_\mathbf{B}(F,n) := \{ X \in \mathcal{M}^{(w)}(F) \, 
        \colon 
        \, \lvert L^\ast(X,n) \rvert \leq B_1, \ \dim(X) \leq B_2 \}
    \end{equation}
    is finite.
    In other words, the pair of functions $\{\lvert L^\ast(\cdot,n) \rvert,\dim \colon \mathcal{M}^{(w)}(F) \to \mathbb{R}_{\geq 0} \}$ has the Northcott property.
\end{theorem}

The proof of \cref{thm:Northcott_left_critical_strip} is divided in two parts. First of all, we show in \cref{prop:Selberg_bound} that certain axiomatic properties of $L$-functions, which are introduced in \cref{def:L_class} and bear some similarity to those used to define the celebrated Selberg class (as we recall in \cref{rmk:Selberg}), are sufficient to guarantee that bounding one of the special values of these $L$-functions which lie at the left of the critical strip yields a bound for the conductor of the $L$-function in question.
Then, we conclude the proof of \cref{thm:Northcott_left_critical_strip} in \cref{sec:Northcott_pure_motives}, by using a result of Deligne \cite{Deligne_1985} on the Northcott property of conductors of Galois representations (see \cref{ex:conductor_l-adic} for further details).

We remark that the set $\mathcal{M}^{(w)}(F)$ should coincide with the set of isomorphism classes of all motives $X \in \mathcal{M}_F$ which are pure of weight $w$, as we explain in \cref{rmk:motivic_L}, but proving this seems out of reach for the mathematics of today. Nevertheless, we show in \cref{cor:modular_abelian_varieties_left} that \cref{thm:Northcott_left_critical_strip} implies an unconditional result for those abelian varieties which are known to be potentially modular, such as elliptic curves and abelian surfaces defined over totally real fields. 

Another case in which the expected properties of $L$-functions are known is provided by the Dedekind $\zeta$-functions $\zeta_F(s)$ associated to number fields $F$. In particular, these can be seen as the $L$-functions associated to the Weil restrictions, from $F$ to $\mathbb{Q}$, of the trivial motive $\mathbf{1}_F \in \mathcal{M}_F$.
Therefore, it is natural to study the Northcott properties for the special values of Dedekind $\zeta$-functions at the integers. We summarise the results that we achieve in the following theorem, which shows in particular that the special values of Dedekind $\zeta$-functions at an integer $n \in \mathbb{Z}$ satisfy the Northcott property if and only if $n \leq 0$. 

\begin{theorem}
\label{thm:Northcott_class_number_formula}
Let $S(\mathbb{Q})$ be the set of isomorphism classes of number fields. For every real number $B \in \mathbb{R}_{\geq 0}$ and every $n \in \mathbb{Z}$, we set $S_{B}(\mathbb{Q},n) := \{ [F] \in S(\mathbb{Q}) \colon \vert\zeta_F^{\ast}(n)\vert \leq B \}$. Then:
\begin{itemize}
    \item
for every $n \in \mathbb{Z}_{\geq 1}$, there exists $B_n \in \mathbb{R}_{\geq 0}$ such that for each $B \geq B_n$ the set $S_B(\mathbb{Q},n)$
is infinite;
\item for every $n \in \mathbb{Z}_{\leq 0}$ and $B \in \mathbb{R}_{> 0}$, the set $S_B(\mathbb{Q},n)$ is finite. 
\end{itemize}

Moreover, there exists an absolute, effectively computable constant $c_0 \in \mathbb{R}_{> 0}$ such that the following upper bound
\begin{equation} \label{eq:intro_explicit_n}
    \lvert S_{B}(\mathbb{Q},n) \rvert \leq \exp\left( \frac{c_0}{1 - n} \log\left( B \right) \left( \log\log\left( B \right) \right)^3 \right)
\end{equation}
holds true for every $n \in \mathbb{Z}_{\leq -1}$ and every $B \in \mathbb{R}_{> 1}$. 
Finally, there exist two other absolute, effectively computable constants $c_1, c_2 \in \mathbb{R}_{> 0}$ such that the following upper bound
\begin{equation} \label{eq:intro_explicit_0} 
    \lvert S_B(\mathbb{Q},0) \rvert \leq \exp\left( c_1 B^{c_2 \log\log(B)} (\log\log(B))^3 \right)
\end{equation}
holds true for every $B \in \mathbb{R}_{> 1}$.
\end{theorem}

We divide the proof of \cref{thm:Northcott_class_number_formula} in various steps. First of all, we devote \cref{prop:Dedekind_left} to the study of the sets $S_B(\mathbb{Q},n)$ for integers $n \in \mathbb{Z} \setminus \{0,1\}$. In particular, despite the fact that we have a natural inclusion $S(\mathbb{Q}) \subseteq \mathcal{M}^{(0)}(\mathbb{Q})$, given by sending $F$ to the motive $\mathbf{1}_{F/\mathbb{Q}}$ obtained as the Weil restriction of the trivial motive $\mathbf{1}_F \in \mathcal{M}_F$, we cannot apply directly \cref{thm:Northcott_left_critical_strip} to prove that the sets $S_B(\mathbb{Q},n)$ are finite when $n \leq -1$. Indeed, the dimension of $\mathbf{1}_{F/\mathbb{Q}}$ equals the degree $d_F := [F \colon \mathbb{Q}]$, which is a priori not bounded when we bound the special values of Dedekind $\zeta$-functions. In order to overcome this issue, we need to adapt the proof of \cref{thm:Northcott_left_critical_strip} to this setting, using the fact that the discriminant $\lvert \Delta_F \rvert$ of a number field $F$ grows exponentially with its degree, as we recall in \cref{ex:discriminants}.
Moreover, the explicit upper bound displayed in \eqref{eq:intro_explicit_n} is obtained by finding an explicit upper bound for the discriminant of the number fields belonging to $S_B(\mathbb{Q},n)$, which is used in combination with an explicit version of Hermite's theorem, proved by Couveignes in \cite{Couveignes_2020} (see also \cref{ex:discriminants}).

To conclude the proof of \cref{thm:Northcott_class_number_formula}, we devote \cref{sec:boundary_NF} to the study of special values of Dedekind $\zeta$-functions at the boundary of the critical strip. In particular, \cref{prop:Northcott_class_number_formula} shows that the Northcott property holds for $n=0$ and does not hold for $n=1$, whereas  \cref{prop:Northcott_class_number_formula_quantitative} is devoted to the proof of the explicit upper bound \eqref{eq:intro_explicit_0}.  
We prove this bound by applying Stark's effective version of the Brauer-Siegel theorem for CM fields \cite{st74} together with Zimmert's explicit lower bounds for the regulator of a number field \cite{zi81}, combined with previous work of the first author of this paper \cite{pa14}. 

Let us conclude this sub-section by noting that \cref{thm:Northcott_left_critical_strip,thm:Northcott_class_number_formula} show how the validity of the Northcott property for the special values of $L$-functions depends crucially on the point at which these special values are taken. 
Such a dependence is also emphasised in recent work of Généreux, Lalín and Li \cite{Genereux_Lalin_Li_2022}, and of Généreux and Lalín \cite{Genereux_Lalin_2022}, which generalises \cref{thm:Northcott_class_number_formula} by looking at the special values at any complex number $s \in \mathbb{C}$ of the Dedekind $\zeta$-functions associated to global fields in any characteristic (see \cref{rmk:Lalin} for more details).  

\subsection{Outlook} To conclude our paper, we devote \cref{sec:abelian_varieties} to $L$-functions associated to abelian varieties, and in particular to the study of their special values at the centre of the critical strip.
More precisely, we show that even for elliptic curves over $\mathbb{Q}$ it is not clear whether or not this special value would satisfy a Northcott property, even if we assume some of the deepest conjectures concerning this number, such as the Birch and Swinnerton-Dyer conjecture.

\section{Properties of heights}
\label{sec:properties_of_heights}
In complete generality we may say that a \textit{height} (or \textit{height function}) on a set $S$ is a function $h \colon S \to \Gamma$ with values in a partially ordered set $\Gamma$.
The aim of this section is to describe various properties of height functions in this generality.

\subsection{Northcott property}
\label{sec:general_definitions_heights}
 The first and more restrictive property of heights we will introduce is named after Northcott's theorem \cite{no49}, which shows the finiteness of sets of points whose height and degree is bounded (see \cref{ex:Weil_height}).
    
    \begin{definition} \label{def:northcott}
        Let $h \colon S \to \Gamma$ be a height function, and let $\mathbb{S}$ be a collection of subsets of $S$.
        Then the height $h$ has:
        \begin{itemize}
            \item[(i)] the \textit{fibre-wise $\mathbb{S}$-Northcott property} if and only if the fibres of $h$ lie in $\mathbb{S}$;
            \item[(ii)] the \textit{$\mathbb{S}$-Northcott property} if and only if $\{s \in S \mid h(s) \leq \gamma\} \in \mathbb{S}$ for every $\gamma \in \Gamma$. 
        \end{itemize}
        When $\mathbb{S}$ is the collection of finite subsets of $S$ it will usually be omitted from the notation.
        \end{definition}
        
        The previous definition readily generalises to sets of height functions, using the following notion of product height. 
        \begin{definition} \label{def:product_height}
            If $\mathbf{h} = \{ h_i \colon S \to \Gamma_i \}_{i \in I}$ is a set of height functions, we define their product $\widetilde{\mathbf{h}}$ as the function:
            \[
                \begin{aligned}
                    \widetilde{\mathbf{h}} \colon S &\to \prod_{i \in I} \Gamma_i \\
                s &\mapsto (h_i(s))_{i \in I}
                \end{aligned}
            \]
            where the set $\prod_{i \in I} \Gamma_i$ is endowed with the product order.
        \end{definition}
        \begin{definition}
        If $\mathbf{h} = \{ h_i \colon S \to \Gamma_i \}_{i \in I}$ is a set of height functions we say that $\mathbf{h}$ has one of the properties described in \cref{def:northcott} if and only if the product height $\widetilde{\mathbf{h}}$  has these properties.
    \end{definition}
    Before moving on, let us observe that
    \begin{equation}
    \label{eq:Northcott_implies_fibre-wise_Northcott}
        \boxed{h \ \text{has} \ \mathbb{S}\text{-Northcott}} + \boxed{\mathbb{S} \ \text{is lower-closed}} \Rightarrow \boxed{h \ \text{has} \ \text{fibre-wise} \ \mathbb{S}\text{-Northcott}}
    \end{equation}
    where $\mathbb{S}$ is called \textit{lower-closed} if for all $Y \subseteq X \subseteq S$ then $X \in \mathbb{S} \Rightarrow Y \in \mathbb{S}$.
    Moreover, if $\mathbb{S}$ is the collection of finite subsets of $S$ then
    \[
        \boxed{h \ \text{has} \ \text{fibre-wise Northcott}} + \boxed{h(S) \ \text{is upper-finite}} \Rightarrow \boxed{h \ \text{has} \ \text{Northcott}} 
    \]
    where we say that $X \subseteq \Gamma$ is \textit{upper-finite} if $X_{\leq \gamma} := \{x \in X \mid x \leq \gamma\}$ is finite for all $\gamma \in \Gamma$.

    \subsection{Bogomolov property}
    \label{sec:Bogomolov}
    Let us now shift to the definition of the Bogomolov property.
    This uses the concepts of \textit{essential infimum} and \textit{successive infima}, that we now review.
    
    \begin{definition} \label{def:essential_minimum}
        Let $\Gamma$ be a partially ordered set, let $X \subseteq \Gamma$ and let $\mathbb{X}$ be a collection of subsets of $X$. 
        Write $X_{\leq \gamma} := \{ x \in X \mid x \leq \gamma \}$ for every $\gamma \in \Gamma$. 
        Then $X$ has an \textit{$\mathbb{X}$-essential infimum} (resp. \textit{$\mathbb{X}$-essential minimum}) if the set
        \[
            \Xi(X,\mathbb{X}) := \{ \gamma \in \Gamma \mid X_{\leq \gamma} \not\in \mathbb{X} \} \subseteq \Gamma
        \]
        has an infimum (resp. a minimum). 
        In this case we denote the set of infima (resp. minima) of $\Xi(X,\mathbb{X})$ by $\mu_{\text{ess}}(X,\mathbb{X}) \subseteq \overline{\Gamma}$. 
        Here $\overline{\Gamma} := \Gamma \sqcup \{+\infty\}$ is the partially ordered set obtained by adjoining to $\Gamma$ a global maximum $+\infty$. In particular, $\mu_{\text{ess}}(X,\mathbb{X}) = \{ +\infty \}$ if and only if $\Xi(X,\mathbb{X}) = \emptyset$, \textit{i.e.} if and only if $X_{\leq \gamma} \in \mathbb{X}$ for every $\gamma \in \Gamma$.
    \end{definition}

    \begin{definition} \label{def:successive_minima_least}
    Let $\Gamma$ be a partially ordered set and let $k \in \mathbb{N}$. Then a subset $X \subseteq \Gamma$ has \textit{at least $k$ successive sets of infima} (respectively \textit{at least $k$ successive sets of minima}) if:
    \begin{itemize}
        \item[(i)] $X$ is bounded from below;
        \item[(ii)] whenever $k \geq 1$, $X$ has at least $k - 1$ successive sets of infima (resp. sets of minima) and the set $X \setminus X_{k - 1}$ has an infimum (resp. minimum). In this case, we denote by $\mu_{k}(X) \subseteq \Gamma$ the set of infima (resp. minima) of $X \setminus X_{k - 1}$. 
        Moreover, the set $X_{k - 1}$ is defined by induction as $X_0 := \emptyset$ and
        \[
            X_{k - 1} := X_{k - 2} \cup U_{k - 1}
        \]
        for any $k \geq 2$, where $U_{k - 1} \subseteq \Gamma$ denotes the union of connected components of the set $X \cup \{\mu_{k - 1}(X)\}$ that contain an element of $\mu_{k - 1}(X)$.
        These connected components are taken with respect to the subspace topology induced on $X \cup \{\mu_{k - 1}(X)\}$ by the order topology on $\Gamma$.
    \end{itemize}
    \end{definition}
    It is easy to see that for every $j \in \mathbb{Z}_{\geq 1}$ and every $x_j \in \mu_j(X)$ and $x_{j + 1} \in \mu_{j + 1}(X)$ the inequality $x_j \leq x_{j + 1}$ holds. 
    Moreover, if $\mu_{j+1}(X) = \mu_{j}(X)$ for some $j \in \mathbb{Z}_{\geq 1}$ then $X$ has at least $k$ successive infima for every $k \in \mathbb{N}$ and $\mu_k(X) = \mu_j(X)$ for every $k \geq j$. This leads to the following definition.
    
    \begin{definition}
    \label{def:successive_minima}
        Let $\Gamma$ be a partially ordered set. Then any subset $X \subseteq \Gamma$ has \textit{exacly $k$ successive sets of infima} (respectively \textit{exacly $k$ successive sets of minima}) for some $k \in \mathbb{N}$ if it has at least $k$ successive sets of infima (resp. sets of minima) and at least one of the following holds:
        \begin{itemize}
            \item[(i)] $X$ does not have at least $k+1$ successive sets of infima (resp. sets of minima);
            \item[(ii)] $\mu_{k + 1}(X) = \mu_k(X)$.
        \end{itemize}
    \end{definition}

    The previous definitions, albeit quite abstract, allow us to recover the usual notions of successive minima, appearing in geometry of numbers and in Arakelov geometry, as we explain in \cref{sec:examples_minima}.
    For now, we will instead use the notion of successive minima to give the definition of Bogomolov property.
    \begin{definition}
        Let $h \colon S \to \Gamma$ be a height function. Then $h$ has \textit{Bogomolov number} $\mathcal{B}(h) \in \mathbb{N}$ if the set $h(S) \subseteq \Gamma$
        has exactly $\mathcal{B}(h)$ successive sets of infima, denoted by $\mu_j(h)$ for every $j \in \{1,\dots,\mathcal{B}(h)\}$.
    \end{definition}
    \begin{definition} \label{def:Bogomolov}
        Let $h \colon S \to \Gamma$ be a height function. Then $h$ has:
        \begin{itemize}
            \item[(i)] the \textit{very weak Bogomolov property} if and only if $\mathcal{B}(h) \geq 0$, \textit{i.e.} if and only if the set $h(S) \subseteq \Gamma$ is bounded from below;
            \item[(ii)] the \textit{weak Bogomolov property} if and only if $\mathcal{B}(h) \geq 1$ and $\mu_1(h) \subseteq h(S)$, \textit{i.e.} if and only if $h(S)$ has at least one minimum;
            \item[(iii)] the \textit{Bogomolov property} if and only if either $\lvert h(S) \rvert = 1$ or $\mathcal{B}(h) \geq 2$ and $\mu_1(h) \subseteq h(S)$, \textit{i.e.} if and only if the minima of $h(S)$ are isolated.
        \end{itemize}
        Moreover, for any collection $\mathbb{S}$ of subsets of $S$ the height $h$ has:
        \begin{itemize}
            \item[(iv)] the $\mathbb{S}$-\textit{essential Bogomolov property} if the set $h(S) \subseteq \Gamma$ has an $h(\mathbb{S})$-essential infimum.
        \end{itemize}
    \end{definition}

\begin{remark}
    We could give another natural definition of Bogomolov property, inspired by the theorems of Ullmo \cite{Ullmo_1998} and Zhang \cite{Zhang_1995} on the Bogomolov conjecture. More precisely, given a height $h \colon S \to \Gamma$, the order topology on $\Gamma$ induces a topology on $S$, which is the coarsest topology making $h$ continuous. Then, we can say that $S$ has the \textit{topological Bogomolov property} if this topology coincides with the discrete topology. 
\end{remark}
    
The previous definition readily generalises to sets of height functions.
\begin{definition}
    If $\mathbf{h} = \{h_i \colon S \to \Gamma_i\}_{i \in I}$ is a set of height functions, we write $\mathcal{B}(\mathbf{h})$ and $\mu_{\text{ess}}(\mathbf{h},\mathbb{S})$ for the Bogomolov number and the essential infimum of the product height $\widetilde{\mathbf{h}}$ introduced in \cref{def:product_height}. Moreover, we say that $\mathbf{h}$ has one of the various Bogomolov properties introduced in \cref{def:Bogomolov} if and only if $\widetilde{\mathbf{h}}$ does. 
\end{definition}

Clearly one has the chains of implications
\begin{align*}
    \boxed{h \ \text{has Bogomolov}} \Rightarrow \boxed{h \ \text{has weak Bogomolov}} &\Rightarrow \boxed{h \ \text{has very weak Bogomolov}} \\
    \boxed{h \ \text{has} \ \text{Northcott}} &\Rightarrow \boxed{h \ \text{has Bogomolov}} \, .
\end{align*}

\subsection{Lehmer property}
As we will see in various examples, sometimes one needs to modify the height under consideration in order to obtain a Bogomolov property. The correct modification for the height of algebraic numbers should consist in multiplying it by the degree, as firstly observed by Lehmer \cite{Lehmer_1933}. Therefore, it seems natural to introduce a similar property for general heights.

\begin{definition}
    Let $\mathbf{h} = \{h_i \colon S \to \Gamma_i\}_{i \in I}$ be a set of heights, and let $\alpha \colon \prod_{i \in I} \Gamma_i \to \Gamma$ be any map of sets, where $\Gamma$ is a partially ordered set. Then the \textit{Lehmer number} $\mathcal{L}(\mathbf{h},\alpha) \in \mathbb{N}$ is defined to be the Bogomolov number of the height
    \[
        S \xrightarrow{\widetilde{\mathbf{h}}} \prod_{i \in I} \Gamma_i \xrightarrow{\alpha} \Gamma
    \]
    and the successive infima of $\alpha(\widetilde{\mathbf{h}}(S))$ are denoted by $\mu_j(\mathbf{h},\alpha)$ for $j \in \{1,\dots,\mathcal{L}(\mathbf{h},\alpha)\}$. Moreover, the pair $(\mathbf{h},\alpha)$ has:
    \begin{itemize}
        \item[(i)] the \textit{very weak Lehmer property} if and only if $\alpha \circ \widetilde{\mathbf{h}}$ has the very weak Bogomolov property;
        \item[(ii)] the \textit{weak Lehmer property} if and only if $\alpha \circ \widetilde{\mathbf{h}}$ has the weak Bogomolov property;
        \item[(iii)] the \textit{Lehmer property} if and only if $\alpha \circ \widetilde{\mathbf{h}}$ has the Bogomolov property.
    \end{itemize}
\end{definition}

\begin{remark}
    Note that Lehmer's question, which inspired our definition of the Lehmer property, concerns very specifically the Weil height of algebraic numbers (see \cref{ex:Weil_height}). 
    Nevertheless, it seems natural to introduce a specific terminology to name those heights which satisfy the Bogomolov property only after a suitable modification.
    Note in particular that this framework encompasses also the analogues of Lehmer's problem which have been considered for the Néron-Tate height of points on an abelian variety, as we recall \cref{ex:canonical_height}.
\end{remark}

It is easy to observe that we have the following implications
\begin{align*}
    \boxed{h' \ \text{has very weak Bogomolov}} + \boxed{\alpha \circ \widetilde{\mathbf{h}} \geq h'} \Rightarrow& \boxed{(\mathbf{h},\alpha) \ \text{has very weak Lehmer}} \\
    \boxed{h' \ \text{has weak Bogomolov}} + \boxed{\alpha \circ \widetilde{\mathbf{h}} \geq h'} \Rightarrow& \boxed{(\mathbf{h},\alpha) \ \text{has weak Lehmer}} \\
    \boxed{h' \ \text{has Bogomolov}} + \boxed{\alpha \circ \widetilde{\mathbf{h}} \geq h'} \Rightarrow& \boxed{(\mathbf{h},\alpha) \ \text{has Lehmer}}
\end{align*}
where $h' \colon S \to \Gamma$ is any height and $\alpha \circ \widetilde{\mathbf{h}} \geq h'$ means that $\alpha( \widetilde{\mathbf{h}}(s)) \geq h'(s)$ for every $s \in S$.

\section{Examples of successive infima}
\label{sec:examples_minima}

We devote this short section to the study of examples of successive infima and minima.
In particular, we will see that our \cref{def:successive_minima_least,def:successive_minima} recover the notions of successive infima and minima present in Arakelov geometry, due to Minkowski (for lattices) and Zhang (for heights associated to hermitian line bundles).

\begin{example}
\label{ex:successive_minima_R}
        Let $\Gamma = \mathbb{R}$. In this case the order topology coincides with the Euclidean topology. Then every set which has at least zero successive infima (\textit{i.e.} is bounded from below) has also at least $n$ successive infima for every $n \in \mathbb{N}$.
        Moreover, if $X \subseteq \mathbb{R}$ is a finite union of open intervals $X = \bigcup_{i = 1}^k (a_i,b_i)$ with $a_1 < b_1 < a_2 < b_2 < \dots$, then it is easy to see that $X$ has exactly $k$ successive infima, with $\mu_i(X) = a_i$ for every $i \in \{1,\dots,k\}$.
        Finally, if $X \subseteq \mathbb{R}$ is countable then $X$ has exactly $k \in \mathbb{Z}_{\geq 1}$ successive minima if and only if 
        there exists a Cauchy sequence $\{x_n\}_{n \in \mathbb{N}} \subseteq X$ such that $\lvert \{ x \in X \mid x \leq x_n, \ \forall n \in \mathbb{N} \} \rvert = k$.
    \end{example}
\begin{example}[Minkowski]
    Let $\Lambda \subseteq \mathbb{R}^n$ be a lattice, and let $g \colon \mathbb{R}^n \to \mathbb{R}_{\geq 0}$ be any distance function (as defined by Cassels in \cite[Chapter~IV]{Cassels_1997}), \textit{i.e.} any continuous function such that $g(t \mathbf{x}) = \lvert t \rvert g(\mathbf{x})$ for all $t \in \mathbb{R}$. 
    Moreover, for every $\lambda \in \Lambda$ we consider the vector space $V_{g,\lambda} := \left\langle \{ x \in \Lambda \mid g(x) \leq g(\lambda) \} \right\rangle_{\mathbb{R}}$.
    Then, the image of the map
    \[
        \begin{aligned}
            \Lambda &\to \mathbb{R}_{\geq 0} \times \mathbb{N} \\
            \lambda &\mapsto \left(g(\lambda), \dim_{\mathbb{R}}(V_{g,\lambda}) \right) 
        \end{aligned}
    \]
    has exactly $n$ successive infima, which are given by the pairs $(\mu_j(\Lambda,g),j)$ for some sequence
        \[
            0 < \mu_1(\Lambda,g) \leq \mu_2(\Lambda,g) \leq \dots \leq \mu_n(\Lambda,g) < +\infty
        \]
    with $\mu_j(\Lambda,g) \in \mathbb{R}_{> 0}$ for every $j \in \{1,\dots,n\}$.
    The numbers $\{\mu_j(\Lambda,g)\}$ are usually called \textit{successive minima} of the function $g$ on the lattice $\Lambda$ (see for instance \cite[Chapter~VIII]{Cassels_1997}). However, these numbers are really infima and not minima in general.
\end{example}

\begin{example}[Zhang]
        Let $\mathcal{X} \to \operatorname{Spec}(\mathbb{Z})$ be an arithmetic variety of dimension $d$, as defined in \cite{Zhang_1995}, and let $\mathsf{Cl}(X)$ be the set of closed sub-schemes of the generic fibre $X := \mathcal{X}_{\mathbb{Q}}$. 
        Fix $\overline{\mathcal{L}}$ to be a relatively semi-ample hermitian line bundle on $\mathcal{X}$ with ample generic fibre, and let $h_{\overline{\mathcal{L}}} \colon X(\overline{\mathbb{Q}}) \to \mathbb{R}$ be the associated height. 
        Then, the image of the map
        \begin{align*}
            \mathsf{Cl}(X) &\to \mathbb{R} \times \mathbb{N} \\
                Y &\mapsto \left( \inf\{ h_{\overline{\mathcal{L}}}(x) \, \mid \, x \in X(\overline{\mathbb{Q}}) \setminus Y(\overline{\mathbb{Q}}) \} \ , \ \dim(Y) \right)
        \end{align*}
        has exactly $d+1$ successive infima, which are given by pairs $(\mu_j(\mathcal{X},\overline{\mathcal{L}}),j)$ for some sequence
        \[
            \mu_0(\mathcal{X},\overline{\mathcal{L}}) \leq \mu_1(\mathcal{X},\overline{\mathcal{L}}) \leq \dots \leq \mu_d(\mathcal{X},\overline{\mathcal{L}}) \leq +\infty
        \]
        with $\mu_j(\mathcal{X},\overline{\mathcal{L}}) \in \mathbb{R}$ for every $j \in \{0,\dots,d - 1\}$ and $\mu_d(\mathcal{X},\overline{\mathcal{L}}) \in \mathbb{R} \sqcup \{+\infty\}$.
        It is easy to see that $\mu_d(\mathcal{X},\overline{\mathcal{L}}) = +\infty$ if and only if $X$ is irreducible. Indeed, this happens if and only if for every $Y \in \mathsf{Cl}(X)$ such that $\dim(Y) = d$, we have that $\{h_{\overline{L}}(x) \colon x \in X(\overline{\mathbb{Q}}) \setminus Y(\overline{\mathbb{Q}})\} = \emptyset$, which happens if and only if $X(\overline{\mathbb{Q}}) = Y(\overline{\mathbb{Q}})$.
        
        Moreover, for every $j \in \{0,\dots,d-1\}$ we have $\mu_{j}(\mathcal{X},\overline{\mathcal{L}}) = e_{d - j}(\overline{\mathcal{L}})$, where $e_1(\overline{\mathcal{L}}) \geq \dots \geq e_d(\overline{\mathcal{L}})$ is the sequence defined in \cite[\S~5]{Zhang_1995}. 
\end{example}

\section{Heights and their relations to special values of \texorpdfstring{$L$}{L}-functions}
\label{sec:examples_heights}

The aim of this section is to provide a roundup of examples of heights satisfying the properties introduced in \cref{sec:properties_of_heights}, and to relate these examples to special values of $L$-functions.
Note that a complete understanding of these examples is not necessary to understand the main results of the present paper. In fact, our aim here is to show the variety of situations in which one finds heights having the properties introduced above.

\subsection{Logarithmic Weil height}
\label{ex:Weil_height}
Let us start with the absolute logarithmic Weil height $h \colon \overline{\mathbb{Q}} \to \mathbb{R}$ (see \cite[Definition~1.5.4]{bg06}), which was the main inspiration to give the general definitions that appear in \cref{sec:properties_of_heights}.

It is immediate to see that $h$ does not have the fibre-wise Northcott property (with respect to the collection of finite subsets of $\overline{\mathbb{Q}}$), for example because $h(\zeta) = 0$ for any root of unity $\zeta \in \overline{\mathbb{Q}}$. 
Hence $h$ does not have the Northcott property.
It is also immediate to see that the same holds for the degree $\operatorname{deg} \colon \overline{\mathbb{Q}} \to \mathbb{Z}_{\geq 1}$, where $\operatorname{deg}(\alpha) := [\mathbb{Q}(\alpha) \colon \mathbb{Q}]$ for every $\alpha \in \overline{\mathbb{Q}}$.
However, Northcott's theorem (see \cite[Theorem~1.6.8]{bg06}) shows that the set $\mathbf{h} = \{h,\operatorname{deg}\}$ has the Northcott property.
    Moreover, it is immediate to see that $h$ has the weak Bogomolov property, because $0 \in \mathbb{R}$ is a minimum for $h(\overline{\mathbb{Q}})$, attained exactly at the roots of unity (see \cite[Theorem~1.5.9]{bg06}).
    However, it is easy to see that this minimum is not isolated, because for example $\lim_{n \to +\infty} h(\sqrt[n]{2}) = 0$. Hence $\mathcal{B}(h) = 1$, and $h$ does not have the Bogomolov property.
    Finally, asking whether the set $\mathbf{h} = \{h,\operatorname{deg}\}$ has the Lehmer property with respect to the function 
    \begin{equation} \label{eq:Lehmer_map_pi}
        \begin{aligned}
            \pi \colon \mathbb{R} \times \mathbb{Z}_{\geq 1} &\to     \mathbb{R} \\
            (x,d) &\mapsto x  d
        \end{aligned}
    \end{equation}
    is equivalent to Lehmer's celebrated problem (see \cite[\S~1.6.15]{bg06}).
    
    We remark that Smyth's theorem \cite[Theorem~4.4.15]{bg06} says that $(\mathbf{h},\pi)$ has the Lehmer property when restricted to the set $S \subseteq \overline{\mathbb{Q}}$ of algebraic numbers which are not Galois-conjugate to their inverse.
    Moreover, Dobrowolski's theorem \cite[Theorem~4.4.1]{bg06} says that, if we let
    \begin{align*}
        \alpha \colon \mathbb{R} \times \mathbb{Z}_{\geq 1} &\to \mathbb{R} \\
        (x,d) &\mapsto x  d  \left( \frac{\log(3 d)}{\log\log(3d)} \right)^3
    \end{align*}
    then the pair $(\mathbf{h},\alpha)$ has Lehmer's property.
    
    Let us mention some of the recent work concerning Northcott, Bogomolov and Lehmer properties relative to the logarithmic Weil height. 
    First of all, it is known that $h$ has the Northcott or Bogomolov property, when restricted to suitable sub-fields of $\overline{\mathbb{Q}}$ having infinite degree over $\mathbb{Q}$.
    We refer the interested reader to \cite{bz01,Dvornicich_Zannier_2008,Widmer_2011,Checcoli_Widmer_2013,Fehm_2018,Checcoli_Fehm_2021} for the study of fields having the Northcott property, and to \cite{bz01,Amoroso_Dvornicich_2000,Habegger_2013,Pottmeyer_2013,adz14,Grizzard_2015,Galateau_2016,Frey_2021,Frey_2022,Plessis_2019} for the study of fields having the Bogomolov one. 
    Moreover, it has also been shown that many fields $F \subseteq \overline{\mathbb{Q}}$ do not have the Bogomolov property. 
    This happens for $F = \mathbb{Q}^\text{tr}(i)$ (as proved by May \cite[Example~1]{May_1972}, and later rediscovered by Amoroso and Nuccio \cite{Amoroso_Nuccio_2007} and by Amoroso, David and Zannier \cite[\S~5]{adz14}), where $\mathbb{Q}^\text{tr}$ denotes the compositum of all the totally real number fields inside $\overline{\mathbb{Q}}$. More generally, each field $F \subseteq \overline{\mathbb{Q}}$ containing an infinite sequence of roots of a number $\alpha \in F$ which is not a root of unity will not satisfy the Bogomolov property.
    It has been proved by Amoroso (see \cite[Theorem~3.3]{Amoroso_2016}) that certain fields $F \subseteq \overline{\mathbb{Q}}$ of this second kind have the Lehmer property, \textit{i.e.} that the pair of functions $\mathbf{h} = \{h,\operatorname{deg}\}$, when restricted to $F$, has the Lehmer property with respect to the function $\pi$ defined in \eqref{eq:Lehmer_map_pi}. 
    This construction has been generalised by Plessis (see \cite[Théorème~1.8]{Plessis_2022}), and is related to a very general conjecture formulated by Rémond in \cite{Remond_2017}, which has been the subject of some recent works, such as \cite{Grizzard_2017} and \cite{Plessis_2023}.
    Finally, recent work of the first author together with Technau and Widmer (see \cite[Theorem~4]{Pazuki_Technau_Widmer_2022}), proves that for every real number $0 < \varepsilon < \gamma \leq 1$ there exists a sub-field $F \subseteq \overline{\mathbb{Q}}$ such that the function $h_\gamma \colon \overline{\mathbb{Q}} \to \mathbb{R}$, defined by $h_\gamma(x) := (\deg(x))^\gamma \, h(x)$, has the Northcott property when restricted to $F$, and for every $\gamma' < \gamma$ the function $h_{\gamma'}$ does not have the Bogomolov property, even when restricted to $F$.
    
    Let us conclude by observing that the work of Akhtari and Vaaler \cite{Akhtari_Vaaler_2016}, combined with the class number formula \eqref{eq:CNF}, show that for every number field $F$ with unit rank $r_F := \operatorname{dim}_{\mathbb{Q}}(\mathcal{O}_F^{\times} \otimes_{\mathbb{Z}} \mathbb{Q})$ there exist units $\{ \gamma_1,\dots,\gamma_{r_F} \} \subseteq \mathcal{O}_F^{\times}$, which form a basis of $\mathcal{O}_F^\times \otimes_\mathbb{Z} \mathbb{Q}$, such that the following inequality
    \begin{equation} \label{eq:zeta_weil_height}
        \frac{h_F \, d_F^{r_F}\,  (2 r_F)!}{2 w_F  (r_F!)^4} \, \prod_{i = 1}^{r_F} h(\gamma_i) \leq \lvert \zeta_F^{\ast}(0) \rvert \leq \frac{h_F \, d_F^{r_F}}{w_F}  \prod_{i = 1}^{r_F} h(\gamma_i)
    \end{equation}
    holds true. This shows that the special value $\zeta_F^{\ast}(0)$ of the Dedekind $\zeta$-function associated to a number field $F$ is commensurable to a product of Weil heights. 
    The terms appearing in \eqref{eq:zeta_weil_height} are given by the degree $d_F := [F \colon \mathbb{Q}]$ of the number field $F$, and by the class number $h_F := \lvert \operatorname{Pic}(\mathcal{O}_F) \rvert$ and the number of roots of unity $w_F := \lvert (\mathcal{O}_F^{\times})_{\text{tors}} \rvert$ of the ring of integers $\mathcal{O}_F$.

\subsection{Mahler measure}
\label{ex:Mahler_measure}

Now, let us study the higher-dimensional generalisation of the function
\begin{align*}
    \pi \circ \widetilde{\mathbf{h}} \colon \overline{\mathbb{Q}} &\to \mathbb{R} \\
    \alpha &\mapsto h(\alpha) \deg(\alpha)
\end{align*}
appearing in \cref{ex:Weil_height}.

To do so, let $\mathbb{G}_{m,\mathbb{C}}^{\infty} := \varprojlim_{n \in \mathbb{N}} \mathbb{G}_{m,\mathbb{C}}^n$ denote the inverse limit of the complex algebraic tori $\mathbb{G}_{m,\mathbb{C}}^n$ with respect to the projections on the last coordinate. Then the global sections of the structure sheaf $\mathcal{O}_{\mathbb{G}_{m,\mathbb{C}}^{\infty}}$ are given by the ring of Laurent polynomials in any number of variables. Moreover, the \textit{logarithmic Mahler measure} is defined as
    \begin{align*}
        m \colon \Gamma(\mathbb{G}_{m,\mathbb{C}}^{\infty},\mathcal{O}_{\mathbb{G}_{m,\mathbb{C}}^{\infty}}) &\to \mathbb{R} \\
        P &\mapsto \int_{\mathbb{T}^\infty} \log\lvert P \rvert \, d\mu_{\mathbb{T}^\infty}.
    \end{align*}
    where $\mathbb{T}^\infty := \varprojlim_{n \in \mathbb{N}} \mathbb{T}^n$ denotes the inverse limit of the real analytic tori $\mathbb{T}^n := (S^1)^n$ with respect to the projections on the last coordinates, and $\mu_{\mathbb{T}^{\infty}}$ denotes the unique Haar probability measure on $\mathbb{T}^{\infty}$.
    
    The height $m$ has the weak Bogomolov property if one restricts it to the ring 
    \[
        \Gamma(\mathbb{G}_{m,\mathbb{Z}}^{\infty},\mathcal{O}_{\mathbb{G}_{m,\mathbb{Z}}^{\infty}}) = \mathbb{Z}[x_1^{\pm 1},x_2^{\pm 1},\dots]
    \]
    of Laurent polynomials with integral coefficients, because for every $P \in \mathbb{Z}[x_1^{\pm 1},x_2^{\pm 1},\dots]$ one has that $m(P) \geq 0$ and $m(P) = 0$ if and only if $P$ is a product of cyclotomic polynomials evaluated at monomials, as proved in \cite{Lawton_1977} (see also \cite{Boyd_1981,Smyth_1981}). 
    In particular, $m(P) = m(\widetilde{P})$ for every $P \in \mathbb{Z}[x_1^{\pm 1},\dots]$, where $\widetilde{P} \in \mathbb{Z}[x_1,\dots]$ denotes the polynomial obtained by clearing out the denominators of $P$.
    Finally, if we let 
    \begin{align*}
        \delta \colon \Gamma(\mathbb{G}_{m,\mathbb{C}}^{\infty},\mathcal{O}_{\mathbb{G}_{m,\mathbb{C}}^{\infty}}) &\to \mathbb{Z}_{\geq 1} \\
        P &\mapsto \sum_{i = 1}^{+\infty} i \, \operatorname{deg}_{x_i}(\widetilde{P}) 
    \end{align*}
    then the pair $(m,\delta)$ has the Northcott property, when restricted to $\Gamma(\mathbb{G}_{m,\mathbb{Z}}^{\infty},\mathcal{O}_{\mathbb{G}_{m,\mathbb{Z}}^{\infty}})$.
    Indeed, this follows from \cite{ma62}, which gives the inequality
    \[
        \operatorname{exp}(m(P)) = \operatorname{exp}(m(\widetilde{P})) \geq 2^{-\sum_{i = 1}^{+\infty} \operatorname{deg}_{x_i}(\widetilde{P})} \, \sum_{\mathbf{j}} \lvert a_{\mathbf{j}} \rvert
    \]
    where $\{a_{\mathbf{j}}\}_{\mathbf{j}} \subseteq \mathbb{Z}$ are the coefficients of $\widetilde{P} = \sum_{\mathbf{j}} a_{\mathbf{j}} x^{a_{\mathbf{j}}}$ written in multi-index notation.
    
    Conjectural relations between the Mahler measure and special values of $L$-functions come from the work of Boyd \cite{bo98} (see also \cite{Brunault_Zudilin_2020} for a survey). 
    For instance, one expects that for every $k \in \mathbb{Z} \setminus \{0,\pm 4\}$, there should exist a rational number $\alpha_k \in \mathbb{Q}^{\times}$ such that
    \begin{equation} \label{eq:deninger_family_Mahler}
        L'(E_k,0) = \alpha_k \, m\left( x + \frac{1}{x} + y + \frac{1}{y} + k \right)
    \end{equation}
    where $E_k$ is the elliptic curve with Weierstrass model $y^2 + k x y = x^3 - 2 x^2 + x$.
    Moreover, it is reasonable to expect that $\alpha_k \in \mathbb{Z}$ for all but finitely many $k$, as suggested by the computational evidence gathered in \cite{bo98}. 
    If this is true, then the relation \eqref{eq:deninger_family_Mahler} and the Northcott property of the Mahler measure would entail a Northcott property for the function $\mathbb{Z} \setminus \{0,\pm 4\} \to \mathbb{R}$ defined by $k \mapsto \lvert L'(E_k,0) \rvert$.
    Note that such a Northcott property holds unconditionally, and follows from \cref{cor:modular_abelian_varieties_left}.

\subsection{Canonical height}
\label{ex:canonical_height}
We have seen that the Mahler measure of a polynomial $P$ can be seen as a way of measuring the complexity of the zero locus of $P$.

There are also more ``canonical'' ways to measure the complexity of sub-varieties of a given arithmetic variety $X$. For instance, given an endomorphism $\phi \colon X \to X$ and a divisor $D \in \operatorname{Div}(X)$ such that $\phi^\ast(D) \sim \alpha D$ for some $\alpha \in \mathbb{R}_{> 1}$, one can construct a ``canonical'' height 
$\widehat{h}^\text{can}_{X,\phi,D} \colon X(\overline{F}) \to \mathbb{R}$, as explained in \cite[\S~B.4]{hs00}. 
In particular, one has a canonical height $\hat{h}_{A,D}$ associated to an abelian variety $A$, polarised by the choice of a divisor $D \in \operatorname{Div}(A)$.

More general notions of canonical heights, which can be applied also to higher dimensional sub-varieties, have been defined by Zhang (see \cite{Zhang_1995,Zhang_1995-Small}), Philippon (see \cite{Philippon_1991,Philippon_1994,Philippon_1995}) and Faltings (see \cite{Faltings_1991} and \cite[Chapter~III, \S~6]{Soule_1992}).
In some cases, these heights turn out to be related to the Mahler measure of some model of the sub-variety in question, as explained for instance in \cite{Maillot_2000} and \cite{Gualdi_2019}.

Let us notice that several Diophantine properties of canonical heights have been studied in the literature. For instance, if $A$ is an abelian variety defined over a number field $F$, and $D \in \operatorname{Div}(A)$ is ample, then it is known that the pair of functions $\{\hat{h}_{A,D}, \deg \colon A(\overline{F}) \to \mathbb{R} \}$ has the Northcott property, as a consequence of Northcott's original work \cite{Northcott_1950}. 
Moreover, as we mentioned already in \cref{sec:Bogomolov}, the work of Ullmo \cite{Ullmo_1998} and Zhang \cite{Zhang_1998} shows that for every sub-variety $X \subseteq A$ the restriction of $\hat{h}_{A,D}$ to $X^\ast(\overline{F})$ has the Bogomolov property, where $X^\ast$ denotes the complement of the union of those positive-dimensional sub-varieties of $X$ which are translates of abelian sub-varieties of $A$ by torsion points. 
Finally, one can formulate an analogue of Lehmer's problem for these heights, as explained in \cite{David_Hindry_2006}, and prove analogues of Dobrowolski's bound in this setting, see for instance \cite{Ratazzi_2004, Carrizosa_2009}. In the analogous case of Drinfeld modules, there is also \cite{Denis_1992, Bosser_Galateau_2019}.

To conclude this subsection, let us recall that canonical heights are related to special values of $L$-functions by a far reaching program initiated by the seminal work of Gross and Zagier (see \cite{Gross_Zagier_1986,Gross_Kohnen_Zagier_1987}, as well as the modern reviews \cite{Conrad_2004,Yuan_Zhang_Zhang_2013}).
    This was later continued by the groundbreaking works of Kudla and collaborators (see \cite{Gross_Kudla_1992,Kudla_1997}). 
    We refer the interested reader to the survey article \cite{Kudla_2004}, as well as to the monograph \cite{Kudla_Rapoport_Yang_2006}.
    Finally, we mention the recent work of Li and Zhang \cite{Li_Zhang_2022_Unitary,Li_Zhang_2022_Orthogonal}, which settles the local Kudla-Rapoport conjecture in the unitary and orthogonal case.

\subsection{Faltings height}
\label{ex:faltings_height}

A particularly interesting example of height of geometric objects has been introduced by Faltings \cite{fa86}, and served as a key tool in his proof of the Mordell conjecture.

More precisely, let $\mathcal{A}(\overline{\mathbb{Q}})$ be the set of isomorphism classes of abelian varieties defined over $\overline{\mathbb{Q}}$, and let $h \colon \mathcal{A}(\overline{\mathbb{Q}}) \to \mathbb{R}$ be the stable Faltings height (see \cite[Section~3]{fa86} and \cite[Page~27]{de85}, which use two different normalizations). Then the set $\{ h, \operatorname{dim} \}$ has the very weak Bogomolov property, since one has the lower bound $h(A) \geq -\log(\sqrt{2 \pi})  \operatorname{dim}(A)$, which is shown by Gaudron and Rémond in \cite[Corollary~8.4]{gr14} following ideas of Bost.
    Moreover, \cite[Page~29]{de85} shows that $h$ has the weak Bogomolov property if we restrict to the set $\mathcal{A}_1(\overline{\mathbb{Q}})$ of $\overline{\mathbb{Q}}$-isomorphism classes of elliptic curves defined over $\overline{\mathbb{Q}}$. Finally, \cite{lo17} and \cite{Burgos_Gil_Menares_Rivera-Letelier_2018} show that $h \colon \mathcal{A}_1(\overline{\mathbb{Q}}) \to \mathbb{R}$ has the Bogomolov property \textit{tout court}. 
    It seems reasonable to ask whether the set $\{ h, \operatorname{dim} \}$ has the Bogomolov property.
    
    Now, let us recall that Faltings's celebrated theorem \cite[Theorem~1]{fa86}, combined with Zarhin's ``trick'' \cite[Remark~16.12]{mi86}, shows that the set $\{ h, \operatorname{dim}, \mathrm{degdef} \}$ has the Northcott property. 
    The ``degree of definition'' function is defined by
    \begin{align*}
        \mathrm{degdef} \colon \mathcal{A}(\overline{\mathbb{Q}}) &\to \mathbb{N} \\
        A &\mapsto \min\{ [F \colon \mathbb{Q}] \mid A \ \text{is defined over $F$} \}
    \end{align*}
    where we say that an abelian variety $A$ defined over a field $\mathbb{L}$ is defined over a sub-field $\mathbb{K}$ if there exists an abelian variety $A'$ defined over $\mathbb{K}$ and such that $A \cong A' \times_{\operatorname{Spec}(\mathbb{K})} \operatorname{Spec}(\mathbb{L})$. Then $\operatorname{degdef}$ is well defined, because every abelian variety defined over $\overline{\mathbb{Q}}$ can be defined over a number field (see \cite[Th{\'e}or{\`e}me~8.8.2]{EGA-IV.3}). 
    It has also been recently proved by Mocz that (if one assumes Artin's and Colmez's conjectures) the function $h$ satisfies Northcott's property, if we restrict to the subset of isomorphism classes of abelian varieties with complex multiplication (see \cite[Theorem~1.4]{mo17}).
    
    For abelian varieties with complex multiplication, the stable Faltings height $h \colon \mathcal{A}(\overline{\mathbb{Q}}) \to \mathbb{R}$ is expected to be related to $L$-functions by Colmez's conjecture \cite[Conjecture~0.4]{co93}, which predicts the relation 
    \begin{equation} \label{eq:Colmez_conjecture}
        - \, h(A) \stackrel{?}{=} \sum_{\chi} m_{(E,\Phi)}(\chi) \, \left( \frac{L'(\chi,0)}{L(\chi,0)} + \log(f_{\chi}) \right)
    \end{equation}
    where $(E,\Phi)$ is the CM-type of $A$ and the sum runs over all the Artin characters
    \[
        \chi \colon G_\mathbb{Q} \to \mathbb{C}
    \] 
    such that $\chi(c) = -1$, where $c \in G_{\mathbb{Q}} := \operatorname{Gal}(\overline{\mathbb{Q}}/\mathbb{Q})$ denotes complex conjugation. This implies in particular that $L(\chi,0) \in \mathbb{C}^{\times}$.
    Moreover, $f_{\chi} \in \mathbb{N}$ denotes the Artin conductor of $\chi$ and the family of rational numbers $\{ m_{(E,\Phi)}(\chi) \}_{\chi} \subseteq \mathbb{Q}$ is defined by the equality
    \[
        \frac{1}{[ G_{\mathbb{Q}} \colon \operatorname{Stab}(\Phi) ]} \sum_{\sigma \in G_{\mathbb{Q}} / \operatorname{Stab}(\Phi)} \lvert \Phi \cap \sigma \circ \Phi \rvert = \sum_{\chi} m_{(E,\Phi)}(\chi) \; \chi(\sigma) 
    \]
    which holds for every $\sigma \in \operatorname{Gal}(\overline{\mathbb{Q}}/\mathbb{Q})$. In particular, $m_{(E,\Phi)}(\chi) = 0$ for all but finitely many Artin characters.

\subsection{Conductors and discriminants}
We have seen in the previous subsection that Colmez's conjectural formula \eqref{eq:Colmez_conjecture} involves the \textit{Artin conductor} $f_\chi \in \mathbb{N}$ associated to an Artin character $\chi \colon G_\mathbb{Q} \to \mathbb{C}$.
By definition $f_\chi := \operatorname{N}(\mathfrak{f}_\rho)$ coincides with the norm of the conductor associated to any complex representation $\rho \colon G_\mathbb{Q} \to \operatorname{GL}_n(\mathbb{C})$ such that $\chi = \operatorname{tr} \circ \rho$, where $\operatorname{tr} \colon \operatorname{GL}_n(\mathbb{C}) \to \mathbb{C}$ denotes the trace map.
The aim of this subsection is to show that the association $\rho \mapsto \mathfrak{f}_\rho$, and more general types of conductors, behave almost like a height.

\begin{example}[Analytic conductor]
    Let $F$ be a number field, fix $n \in \mathbb{N}$ and let $\mathcal{A}_n(F)$ be the set of cuspidal automorphic representations of $\operatorname{GL}_n(\mathbb{A}_F)$ (see \cite[\S~1]{Iwaniec_Sarnak_2010}). Then Brumley has shown in \cite[Corollary~9]{Brumley_2006} that the \textit{analytic conductor} $\mathcal{C} \colon \mathcal{A}_n(F) \to \mathbb{R}_{\geq 1}$, which is defined in \cite[Equation~(31)]{Iwaniec_Sarnak_2010}, satisfies the Northcott property. In particular, the $n = 1$ case shows that the set of Hecke characters $\psi \colon \mathbb{A}_F^{\times} \to \mathbb{C}^{\times}$ with bounded analytic conductor is finite.
\end{example}

\begin{example}[Conductors of complex representations] \label{ex:conductor_complex}
    Let $F$ be a number field, and $\mathcal{W}_{\mathbb{C}}(F)$ be the set of isomorphism classes of pairs $(V,\rho)$ where $V$ is a finite dimensional complex vector space and $\rho \colon W_F \to \operatorname{GL}(V)$ is a continuous representation of the Weil group $W_F$ (see \cite[\S~1]{Tate_1979}) which is \textit{semi-simple}, \textit{i.e.} a direct sum of irreducible representations. Then, there is a function $\mathfrak{f} \colon \mathcal{W}_{\mathbb{C}}(F) \to \mathrm{Div}(\mathcal{O}_F)$, sending each $(V,\rho)$ to its global Artin conductor ideal $\mathfrak{f}_{\rho} \subseteq \mathcal{O}_F$ (see \cite[Chapter~VII, \S~11]{ne99}). 
    Moreover, the Archimedean local Langlands correspondence, explained for example in \cite{Knapp_1994}, allows one to
    associate to each $(V,\rho) \in \mathcal{W}_{\mathbb{C}}(F)$ an Archimedean conductor $\mathcal{C}_{\infty}((V,\rho)) \in \mathbb{R}$, defined in exactly the same way as the Archimedean part of the analytic conductor of a cuspidal automorphic form.
    Then \cite[Theorem~3.3]{abcz94} can be combined with our previous discussion to show that the function $\mathcal{C} \colon \mathcal{W}_{\mathbb{C}}(F) \to \mathbb{R}$ defined as $\mathcal{C}((V,\rho)) := \operatorname{N}_{F/\mathbb{Q}}(\mathfrak{f}_{\rho}) \, \mathcal{C}_{\infty}((V,\rho))$ satisfies the Northcott property. Let us observe that:
    \begin{itemize}
        \item[(i)] one can consider all the number fields at once as follows: if $\mathcal{W}_{\mathbb{C}}$ denotes the set of isomorphism classes of triples $(F,V,\rho)$, where $F$ is a number field and $(V,\rho) \in \mathcal{W}_{\mathbb{C}}(F)$, then \cite[Property~(a2)]{Rohrlich_1994} shows that the composite map $\mathcal{C} \circ \operatorname{Ind} \colon \mathcal{W}_{\mathbb{C}} \to \mathcal{W}_{\mathbb{C}}(\mathbb{Q}) \to \mathbb{R}$ satisfies the Northcott property, where $\operatorname{Ind} \colon \mathcal{W}_{\mathbb{C}} \to \mathcal{W}_{\mathbb{C}}(\mathbb{Q})$ sends $(F,V,\rho)$ to the induced representation on $W_{\mathbb{Q}} \supseteq W_F$;
        \item[(ii)] the conductor $\mathfrak{f}_\rho$ is related to $L$-functions by means of the functional equation (as explained by Tate in \cite[Theorem~3.5.3]{Tate_1979}).
    \end{itemize}
\end{example}

The following example is the analogue of the previous one for representations valued in vector spaces defined over $\mathbb{Q}_\ell$. It will play a crucial role in our proof of \cref{thm:Northcott_left_critical_strip}.

\begin{example}[Conductors of $\ell$-adic representations] \label{ex:conductor_l-adic}
    Let $F$ be a number field, and $M_F^0$ denote its set of finite places.
    For every $w \in \mathbb{N}$ and every prime number $\ell \in \mathbb{N}$, we let $\mathcal{G}_\ell^{(w)}(F)$ be the set of isomorphism classes of pairs $(V,\rho)$ where $V$ is a finite dimensional vector space over $\mathbb{Q}_\ell$ and $\rho \colon G_F := \operatorname{Gal}(\overline{F}/F) \to \operatorname{GL}(V)$ is a continuous representation satisfying the following properties:
    \begin{itemize}
        \item $\rho$ is semi-simple, \textit{i.e.} a direct sum of irreducible representations;
        \item the set $S_\rho \subseteq M_F^0$ of non-Archimedean places at which $\rho$ is ramified is finite;
        \item $\rho$ is pure of weight $w$. In other words, for every $v \in M_F^0 \setminus S_\rho$ the characteristic polynomial of $\rho(\mathrm{Frob}_v)$ has integer coefficients, and each of its roots has absolute value $\lvert \kappa_v \rvert^{w/2}$, where $\kappa_v$ denotes the residue field of the local field $K_v$, and $\mathrm{Frob}_v \subseteq G_F$ denotes the set of geometric Frobenius elements.
    \end{itemize}

    Now, we define two functions $\dim, \mathfrak{M} \colon \mathcal{G}_\ell^{(w)}(F) \to \mathbb{N}$ by setting $\dim(V,\rho) := \dim(V)$ and $\mathfrak{M}(V,\rho) := \max(\{\mathrm{char}(\kappa_v) \colon v \in S_\rho\})$.  Then, \cite[Corollaire~1]{Deligne_1985} shows that the set $\{\dim,\mathfrak{M}\}$ has the Northcott property. 
    In particular, if we let $\mathfrak{f}_\rho \subseteq \mathcal{O}_F$ denote the conductor ideal of a representation $(V,\rho) \in \mathcal{G}_\ell^{(w)}(F)$ (see for instance Ulmer's work \cite{ul16}), and we set $\mathcal{C}_0(V,\rho) := \lvert \operatorname{N}_{F/\mathbb{Q}}(\mathfrak{f}_\rho)\rvert \in \mathbb{N}$, then we see that $\mathfrak{M}(V,\rho) \leq \mathcal{C}_0(V,\rho)$ for every $(V,\rho) \in \mathcal{G}_\ell^{(w)}(F)$. Therefore, the pair of functions:
    \begin{equation*}
            \left\{ \dim, \mathcal{C}_0 \colon \mathcal{G}_\ell^{(w)}(F) \to \mathbb{N} \right\}
    \end{equation*}
    has the Northcott property.
    
    Let us conclude by making the following observations:
    \begin{enumerate}[label=(\alph*)]
        \item the semi-simplifications of the $\ell$-adic \'etale cohomology groups $H^i_{\text{\'et}}(X_{\overline{F}};\mathbb{Q}_\ell(j))$ associated to a smooth and proper variety $X$ defined over $F$ give rise to elements of $\mathcal{G}_\ell^{(i-2j)}(F)$. 
        For these Galois representations the set $S_\rho$ is contained in the set of primes of $F$ which either lie above $\ell$ or are primes of bad reduction for $X$. 
        This follows from the smooth and proper base change theorem for \'etale cohomology, combined with Deligne's proof of the Weil conjectures (as explained by Jannsen in \cite[Appendix~C]{ja90}).
        \item we can consider all the number fields at once, as we did in \cref{ex:conductor_complex}, by defining $\mathcal{G}_\ell^{(w)}$ as the set of isomorphism classes of triples $(F,V,\rho)$ where $F$ is a number field and $(V,\rho) \in \mathcal{G}_\ell^{(w)}(F)$. Then \cite[Property~(a2)]{Rohrlich_1994} implies that the set $\{\dim \, \circ \, \operatorname{Ind},\mathcal{C}_0 \, \circ \, \operatorname{Ind} \}$ has the Northcott property. 
        Here $\operatorname{Ind} \colon \mathcal{G}_\ell^{(w)} \to \mathcal{G}_\ell^{(w)}(\mathbb{Q})$ is again the map sending $(F,V,\rho)$ to the representation induced on $\operatorname{Gal}(\overline{\mathbb{Q}}/\mathbb{Q}) \supseteq \operatorname{Gal}(\overline{F}/F)$;
        \item the conductor $\mathfrak{f}_\rho$ is supposed to be related to the $L$-function $L(\rho,s)$ by means of the conjectural functional equation, as we recall in \cref{sec:Northcott_pure_motives}.
    \end{enumerate}
\end{example}

To conclude this sub-section, we focus our attention on discriminants of number fields, which are classically known to satisfy a Northcott property, thanks to the work of Hermite and Minkowski. Moreover, given a number field $F$, its discriminant $\Delta_F$ can be seen as the conductor of the Galois representation of $G_\mathbb{Q} := \operatorname{Gal}(\overline{\mathbb{Q}}/\mathbb{Q})$ induced by the trivial representation of $G_F := \operatorname{Gal}(\overline{F}/F)$.

\begin{example}[Discriminants of number fields] \label{ex:discriminants}
    Let $S$ denote the set of isomorphism classes of number fields. Then, it is well known that the function:
    \[
    \begin{aligned}
        S &\to \mathbb{R} \\
        [F] &\mapsto \lvert \Delta_F \rvert
    \end{aligned}
    \]
    has the Northcott property, thanks to Hermite's theorem (see \cite[Chapter~III, Theorem~2.16]{ne99}). Moreover, Malle's conjecture on the distribution of number fields with bounded discriminant and given Galois group, which has been slightly corrected in \cite{Turkelli_2015}, implies that for every $d \in \mathbb{Z}_{\geq 1}$ there should exist a constant $\delta_d \in \mathbb{R}_{> 0}$ such that the following asymptotic:
    \begin{equation} \label{eq:Malle_1}
        \lvert \{ [F] \in S \colon \lvert \Delta_F \rvert \leq X, \ [F \colon \mathbb{Q}] = d \} \rvert \overset{?}{\sim} \delta_d \cdot X
    \end{equation}
    holds true when $X \to +\infty$. This clearly implies that the following asymptotic:
    \begin{equation} \label{eq:Malle_2}
        \lvert \{ [F] \in S \colon \lvert \Delta_F \rvert \leq X, \ [F \colon \mathbb{Q}] \leq d \} \rvert \overset{?}{\sim} \left( \sum_{j=1}^d \delta_j \right) \cdot X
    \end{equation}
    should also hold true when $X \to +\infty$. Since we know that bounding $\lvert \Delta_F \rvert$ yields explicit bounds for the degree $d_F$ of a number field (see for instance \cite{Poitou_1977}), the validity of \eqref{eq:Malle_2} should yield also an estimate for the quantity $\lvert \{[F] \in S \colon \lvert \Delta_F \rvert \leq X \} \rvert$, as $X \to +\infty$. However, there are no explicit conjectures for what the precise asymptotic should be (see \cite{MO_315345} for a discussion), which is related to the fact that it is difficult to give explicit expressions for the constants $\delta_d$ appearing in \eqref{eq:Malle_1}. In particular, we refer the interested reader to Bhargava's work \cite{Bhargava_2007} for a conjecture concerning the constants that should appear when counting number fields with maximal Galois group. We note moreover that the asymptotic portrayed in \eqref{eq:Malle_1} is only known when $d \leq 5$, thanks to the work of Davenport-Heilbronn and Bhargava (see \cite{Bhargava_2006} for a survey).

    Nevertheless, Schmidt \cite{Schmidt_1995}, Ellenberg-Venkatesh \cite{Ellenberg_Venkatesh_2006} and Couveignes \cite{Couveignes_2020} have obtained some explicit upper bounds for the cardinality of the set appearing on the left hand side of \eqref{eq:Malle_2}. 
    In particular, \cite[Theorem~2]{Couveignes_2020} asserts that there exists a positive constant $Q \in \mathbb{R}_{\geq 0}$ such that for every $d \geq Q$ the following upper bound:
    \begin{equation} \label{eq:couveignes}
        \lvert \{ [F] \in S \colon \lvert \Delta_F \rvert \leq X, \ [F \colon \mathbb{Q}] = d \} \leq d^{Q d \log^3(d)} X^{Q \log^3(d)}
    \end{equation}
    holds true. 
    One can then adapt the proof of \cite[Theorem~A.1]{Belolipetsky_2007} to show that there exists a positive constant $\mathfrak{C} \in \mathbb{R}_{> 0}$ such that the following upper bound:
    \begin{equation} \label{eq:effective_Hermite}
    \lvert \{ [F] \in S \colon \lvert \Delta_F \rvert \leq X \} \rvert \leq \exp(\mathfrak{C} \log(X) \log\log(X)^3)
    \end{equation}
    holds true for every $X \geq 1$. Note finally that the recent work of Lemke Oliver and Thorne \cite{Lemke-Oliver_Thorne_2022} refines \eqref{eq:couveignes}, by proving that there exists a constant $Q' \in \mathbb{R}_{> 0}$ such that for every $d \in \mathbb{Z}_{\geq 6}$ there exists a constant $\mathfrak{d}_d \in \mathbb{R}_{> 0}$ with the property that the following upper bound:
    \[
        \lvert \{ [F] \in S \colon \lvert \Delta_F \rvert \leq X, \ [F \colon \mathbb{Q}] = d \} \leq \mathfrak{d}_d X^{Q' \log^2(d)}
    \]
    holds true. Such an upper bound could yield to a refinement of \eqref{eq:effective_Hermite}, where one replaces $(\log\log(X))^3$ with $(\log\log(X))^2$,
    if one manages to prove that $\mathfrak{d}_d \leq d^{Q' d \log^2(d)}$ for every $d \geq 6$.
\end{example}

\subsection{Heights coming from geometry}
We conclude this roundup of examples by talking about two more geometric examples of height: the volume of hyperbolic manifolds and the heights of mixed motives defined by Kato.

\begin{example}[Volumes of hyperbolic manifolds] \label{ex:volume_hyperbolic_manifolds}
    Let $\mathcal{H}$ be the set of isomorphism classes of hyperbolic manifolds of finite volume. Then it is conjectured that the volume $\operatorname{vol} \colon \mathcal{H} \to \mathbb{R}_{\geq 0}$ has the Bogomolov property, and that the minimum is attained at an arithmetic hyperbolic manifold $M \cong \mathfrak{h}_n/\Gamma$, where $\Gamma$ is an arithmetic subgroup of the isometry group of the hyperbolic space $\mathfrak{h}_n$ (as explained by Belolipetsky and Emery in \cite{be14}). Then, if we restrict to the set $\mathcal{H}^{\text{ar}}$ of isomorphism classes of arithmetic hyperbolic manifolds, it is conjectured that the set $\mathbf{h} = \{\operatorname{vol},\operatorname{dim},\operatorname{deg}\}$ has the Northcott property, where the degree is defined by $\operatorname{deg}(M) := [\mathbb{Q}(\operatorname{tr}(\pi_1(M)^{(2)})) \colon \mathbb{Q}]$. Here we denote by $\pi_1(M)^{(2)}$ the sub-group generated by the squares, and by $\operatorname{tr} \colon \pi_1(M) \to \mathbb{C}$ the trace map induced from the embedding of $\pi_1(M)$ into the isomorphism group of $\mathfrak{h}_n$.
    This Northcott property has been proved for three dimensional arithmetic hyperbolic manifolds by Jeon \cite{je14}. 
    
    The relations of hyperbolic volumes with special values of $L$-functions comes for example from the formula
    \[
        \zeta_F^{\ast}(-1) \sim_{\mathbb{Q}^{\times}} \operatorname{vol}\left( \frac{\mathfrak{h}_3^{r_2(F)}}{\Gamma} \right) 
    \]
    which holds for any number field $F$. Here $\Gamma$ is a finite-index and torsion-free subgroup of the group $\mathcal{O}^{(1)} \subseteq \mathcal{O}$ of units having norm one in some order $\mathcal{O} \subseteq \mathcal{B}$ in a totally definite quaternion algebra $\mathcal{B} \neq \operatorname{Mat}_{2 \times 2}(K)$ defined over $K$ (as explained by Vignéras in \cite[Example~IV.1.5]{vi80}).
\end{example}

\begin{example}[Heights of motives] \label{ex:Heights_for_motives}
    Let $F$ be a number field, and let $\mathcal{MM}_F$ denote the category of mixed motives defined by Jannsen (see \cite[\S~4]{ja90}).
    Then Kato constructs in \cite{ka18} a series of height functions which measure the complexity of an object $X \in \mathcal{MM}_F$, using the $v$-adic Hodge theory corresponding to any place $v$ of $F$. 
    One of the richest examples of such a height is given by the function
    \begin{equation} \label{eq:Kato_height}
        \begin{aligned}
            h_{\ast,\diamondsuit} \colon \mathcal{MM}_F &\to \mathbb{R} \\
            X &\mapsto  h_\diamondsuit(X) + \sum_{w \in \mathbb{Z}} h_\ast(\operatorname{gr}^\mathcal{W}_w(X)) 
        \end{aligned}
    \end{equation}
    which is the logarithmic version of the height $H_{\ast,\diamondsuit}$ defined in \cite[\S~1.7.1]{ka18}. 
    Here 
    \begin{equation}
        \operatorname{gr}^\mathcal{W}_w(X) := \frac{\mathcal{W}_w(X)}{\mathcal{W}_{w - 1}(X)}    
    \end{equation}
    denotes the graded piece of $X$ with respect to the ascending weight filtration $\mathcal{W}$, and the various heights $h_\ast(\operatorname{gr}^\mathcal{W}_w(X))$ appearing in \eqref{eq:Kato_height} are a generalisation of Faltings's height (see \cref{ex:faltings_height}) to pure motives.
    On the other hand, the height $h_\diamondsuit(X)$ measures the distance between $X$ and the semi-simplification
    \begin{equation} \label{eq:semi_simplification}
        X^\text{ss} := \bigoplus_{w \in \mathbb{Z}} \operatorname{gr}^\mathcal{W}_w(X)
    \end{equation}
    and thus can be seen as a measure of the mixed nature of $X$.
    
    It is extremely interesting to study the Northcott property for the height $h_{\ast,\diamondsuit}$, in view of the many consequences that this would have, which are investigated in \cite[\S~2]{ka18}. In particular, \cite[Proposition~2.1.17]{ka18} shows that the Northcott property for $h_{\ast,\diamondsuit}$ implies the finite generation of motivic cohomology, which would be a motivic analogue of the Mordell-Weil theorem.
    Special instances of the Northcott property for the height $h_{\ast,\diamondsuit}$ have been recently proved to hold by Koshikawa (see \cite{Kosh15,Kosh16}) and Nguyen (see \cite{Nguyen_2022}).
    In particular, Koshikawa shows that $h_{\ast,\diamondsuit}$ has the Northcott property when restricted to the set of pure motives $X$ which are isogenous to a fixed pure motive $X_0$.
    We note that in this case $h_{\ast,\diamondsuit}(X) = h_\ast(X)$ because $h_{\diamondsuit}(X) = 0$ for pure motives $X$.
    Koshikawa's result is reminiscent of the similar Northcott property for the Faltings height (see \cite[\S~4]{fa86}), which allowed Faltings himself to prove the Tate conjecture for abelian varieties.
\end{example}

\section{Special values outside the critical strip}
\label{sec:left_critical_strip}

The aim of this section is to study Northcott properties for the special values of the $L$-functions attached to pure motives. In order to do so, we will first introduce an axiomatic class of $L$-functions, closely related to the Selberg class (see \cref{def:L_class} and \cref{rmk:Selberg}), whose elements satisfy a precise functional equation and have a well defined associated conductor. Then, we will show in \cref{prop:Selberg_bound} that bounding the special values of these $L$-functions on the left of the critical strip yields a bound for the conductor. This suggests that there should be only finitely many such $L$-functions (see \cref{rmk:Northcott_Selberg}), using a suitable Northcott property for the conductor. Such a Northcott property is known to hold for $L$-functions coming from pure motives of bounded dimension, thanks to the work of Deligne \cite{de85} (see also \cref{ex:conductor_l-adic}), and will allow us to prove \cref{thm:Northcott_left_critical_strip} in \cref{sec:Northcott_pure_motives}.
Therefore, this yields unconditional results for all the pure motives whose $L$-functions are known to satisfy the conjectural properties outlined in \cref{def:L_class}, and in particular for motives associated to potentially modular abelian varieties (see \cref{cor:abelian_varieties_outside,cor:modular_abelian_varieties_left}).
We will finally adapt the aforementioned proof of \cref{thm:Northcott_left_critical_strip} to show that the special values of Dedekind $\zeta$-functions at the left of the critical strip satisfy a Northcott property. This is explained in \cref{prop:Dedekind_left}, where we also remark that the special values of Dedekind $\zeta$-functions taken at the right of the critical strip do not satisfy the Northcott property.

\subsection{Bounding the conductor}

The aim of this section is to show how bounding a special value at the left of the critical strip entails a bound for the conductor of a functional equation. This observation depends only on a few axiomatic properties of $L$-functions, which we recall in the following definition.

\begin{definition} \label{def:L_class}
    Let $F$ be a number field, and $w \in \mathbb{Z}$ a fixed integer, called the \textit{weight}. We define a class $\mathcal{L}^{(w)}(F)$ of meromorphic functions $L(s) \colon \mathbb{C} \dashrightarrow \mathbb{C}$ such that:
    \begin{itemize}
        \item there exists a sequence of polynomials $\{ P_\mathfrak{p}(t) \}_\mathfrak{p} \subseteq \mathbb{Q}[t]$, where $\mathfrak{p}$ runs over the non trivial prime ideals of $\mathcal{O}_F$, such that for every $s \in \mathbb{C}$ with $\Re(s) > \frac{w}{2} + 1$, the function $L(s)$ admits the following absolutely convergent Euler product:
        \begin{equation} \label{eq:Euler_product_Selberg}
            L(s) = \prod_{\mathfrak{p}} \frac{1}{P_\mathfrak{p}(\lvert \mathcal{O}_F/\mathfrak{p}\rvert^{-s})};
        \end{equation}
        \item for every non trivial prime ideal $\mathfrak{p} \subseteq \mathcal{O}_F$, let $d$ be the degree of $P_\mathfrak{p}(t)$ and $\{\alpha_{1,\mathfrak{p}},\dots,\alpha_{d,\mathfrak{p}}\}$ denote the roots of $P_\mathfrak{p}$. Then, we assume that $\lvert \alpha_{j,\mathfrak{p}} \rvert \leq \lvert \mathcal{O}_F/\mathfrak{p} \rvert^{w/2}$ for every $j \in \{1,\dots,d\}$;
        \item there exists an integer $N \geq 1$, a sign $\epsilon \in \{\pm 1\}$, a pair of integers $d_1, d_2 \geq 0$ and two tuples of positive integers $\{\mu_j\}_{j=1}^{d_1}$ and $\{\nu_k\}_{k=1}^{d_2}$ with $1 \leq \mu_j, \nu_k \leq w/2$, such that, if we let:
        \[
            L_\infty(s) := \prod_{j=1}^{d_1} \Gamma_{\mathbb{R}}(s - \mu_j) \cdot \prod_{k=1}^{d_2} \Gamma_\mathbb{C}(s - \nu_k)
        \]
        and $\widehat{L}(s) := L(s) \cdot L_\infty(s)$, then we have the following functional equation: 
        \begin{equation} \label{eq:functional_equation}
            \widehat{L}(s) = \epsilon \cdot N^{\frac{w+1}{2} - s} \cdot \widehat{L}(w + 1 - s)
        \end{equation}
        for every $s \in \mathbb{C}$.
        Here, we recall that the functions $\Gamma_\mathbb{R}(s)$ and $\Gamma_\mathbb{C}(s)$ are defined as:
        \begin{equation}
            \begin{aligned}
                \Gamma_\mathbb{R}(s) &:= \pi^{-s/2} \Gamma(s/2) \\
                \Gamma_\mathbb{C}(s) &:= 2 (2 \pi)^{-s} \Gamma(s)
            \end{aligned}
        \end{equation}
        for every $s \in \mathbb{C}$.
    \end{itemize}

Moreover, if $L(s) \in \mathcal{L}^{(w)}(F)$, the numbers $N_L := N$ and $d_L := d_1 + 2 d_2$ are uniquely determined. 
We call $N_L$ the \textit{conductor} of $L(s)$, and $d_L$ its \textit{degree}.    
\end{definition}
\begin{remark} \label{rmk:Selberg}
    The axioms outlined in \cref{def:L_class} are quite similar to the ones which define the Selberg class (surveyed by Perelli in \cite{Perelli_2005}), and also to the ones outlined in more recent work of Farmer, Pitale, Ryan and Schmidt \cite{Farmer_Pitale_Ryan_Schmidt_2019}. 
    The main difference between our axioms and theirs is that we use the arithmetic normalisation for the functional equation, which depends on a fixed weight $w \in \mathbb{Z}$. This weight is also featured in the second axiom, which can be seen as a version of Ramanujan's conjecture for the $L$-function in question. 
\end{remark}

Now, we are ready to see how any bound for $L$-values at the left of the critical strip yields a bound for the conductor, thanks to the axioms outlined in \cref{def:L_class}.

\begin{proposition} \label{prop:Selberg_bound}
    Let $F$ be a number field. Fix two integers $w, n \in \mathbb{Z}$ such that $2 n < w$. Fix moreover two real numbers $B_1, B_2 \in \mathbb{R}_{\geq 0}$. Then, there exists a real number $B_3 \in \mathbb{R}_{\geq 0}$ such that for every $L \in \mathcal{L}^{(w)}(F)$ such that $d_L \leq B_1$ and $L^\ast(n) \leq B_2$, we have that $N_L \leq B_3$.
\end{proposition}
\begin{proof}
    Combining the assumption $L^\ast(n) \leq B_2$ with the functional equation \eqref{eq:functional_equation}, we see that the following bound:
    \begin{equation} \label{eq:bound_functional_equation}
        N_L^{\frac{w+1}{2} - n} \leq \frac{1}{\lvert L^\ast(w+1-n) \rvert} \cdot \left\lvert \frac{L_\infty^\ast(n)}{L_\infty^\ast(w+1-n)} \right\rvert \cdot B_2
    \end{equation}
    holds true. 
    Now, note that $L^\ast(w+1-n) = L(w+1-n) = \prod_\mathfrak{p} P_\mathfrak{p}(\operatorname{N}_{F/\mathbb{Q}}(\mathfrak{p})^{n - w - 1})^{-1}$, thanks to the assumption that $2 n < w$.
    Hence, we see that:
    \begin{equation} \label{eq:bound_right_inverse}
        \lvert L^\ast(w+1-n) \rvert^{-1} = \prod_\mathfrak{p} \prod_{j=1}^{\deg(P_\mathfrak{p})} \lvert \alpha_{j,\mathfrak{p}} - \lvert \mathcal{O}_F/\mathfrak{p}\rvert^{n-w-1} \rvert \leq B_4
    \end{equation}
    for some constant $B_4$ which depends on $n, w$ and $B_1$ (or, to be more precise, on $n, w$ and $d_L$). 

    Moreover, if $m \in \mathbb{Z}$, we have that:
    \[
        \Gamma_{\mathbb{R}}^\ast(m) = \begin{cases}
            \frac{(-4 \pi)^{-k} (-k)!}{(-2 k)!} \ &\text{if} \ m = 2 k + 1 \leq -1 \\
            \frac{(-\pi)^{-k}}{(-k)!} \ &\text{if} \ m = 2 k \leq 0 \\
            \frac{(k-1)!}{\pi^k}, \ &\text{if} \ m = 2 k \geq 2 \\
            \frac{(2k)!}{(4 \pi)^k \, k!}, \ &\text{if} \ m = 2 k + 1 \geq 1
        \end{cases}
    \]
    which implies that there exists a positive constant $B_5$, depending only on $w$ and $n$, such that:
    \[
        \left\lvert \frac{\Gamma_\mathbb{R}^\ast(n - \mu_j)}{\Gamma_\mathbb{R}^\ast(w+1-n-\mu_j)} \right\rvert \leq B_5
    \]
    for every $j \in \{1,\dots,d_1\}$. Analogously, we have that:
    \[
    \Gamma_{\mathbb{C}}^\ast(m) = \begin{cases}
            2 \frac{(-2 \pi)^{-m}}{(-m)!}, \ &\text{if} \ m \leq 0 \\
            2 \frac{(m-1)!}{(2 \pi)^{m}}, \ &\text{if} \ m \geq 1
        \end{cases}
    \]
    which implies that there exists a positive constant $B_6$, depending only on $w$ and $n$, such that:
    \[
        \left\lvert \frac{\Gamma_\mathbb{C}^\ast(n - \nu_k)}{\Gamma_\mathbb{C}^\ast(w+1-n-\mu_k)} \right\rvert \leq B_6
    \]
    for every $k \in \{1,\dots,d_2\}$. Hence, we see that:
    \begin{equation} \label{eq:bound_gamma}
        \left\lvert \frac{L_\infty^\ast(n)}{L_\infty^\ast(w+1-n)} \right\rvert \leq B_7
    \end{equation}
    where $B_7$ is a positive constant which depends only on $w$, $n$ and $d_L$. Hence, we can combine \eqref{eq:bound_right_inverse} and \eqref{eq:bound_gamma} to see that:
    \[
        N_L^{\frac{w+1}{2} - n} \leq B_8
    \]
    where $B_8$ depends only on $B_1$, $B_2$, $w$ and $n$. Therefore, we can take $B_3 := B_8^{\frac{2}{w+1-2n}}$.
\end{proof}

\begin{remark} \label{rmk:Northcott_Selberg}
    We would like to use \cref{prop:Selberg_bound} to prove that for every $w,n \in \mathbb{Z}$ such that $2 n < w$, and every pair of real numbers $B_1, B_2 \in \mathbb{R}_{\geq 0}$, the set: 
\begin{equation} \label{eq:first_Selberg_Set}
    \{L \in \mathcal{L}^{(w)}(F) \colon d_L \leq B_1, \lvert L^\ast(n) \rvert \leq B_2 \}
\end{equation}
is finite. More precisely, \cref{prop:Selberg_bound} implies that the finiteness of the following set:
\begin{equation} \label{eq:second_Selberg_set}
    \{L \in \mathcal{L}^{(w)}(F) \colon d_L \leq B_1, N_L \leq B_2 \}
\end{equation}
is equivalent to the finiteness of \eqref{eq:first_Selberg_Set}. Such a finiteness result is surely to be expected. Indeed: 
\begin{itemize}
    \item the vector space of functions which satisfy a functional equation similar to \eqref{eq:functional_equation} is known to be finite dimensional when $d_L = 1$, thanks to the work of Bochner (see \cite{Bochner_1958} and \cite[Theorem~3.2]{Conrey_Ghosh_1993}). We refer the interested reader to the recent work of Dixit \cite{Dixit_2020}, and to Perelli's survey \cite{Perelli_2017}, for further extensions of Bochner's result.
    \item in the aforementioned vector space, there should be only finitely many elements which satisfy an Euler product. This property is analogous to the fact that in the vector space of cusp forms there are only finitely many Hecke eigenforms. 
\end{itemize}
\end{remark}

In the current paper, we will not address further the finiteness of the set \eqref{eq:second_Selberg_set}. Instead, we will restrict our attention to $L$-functions of ``motivic'' origin, for which one knows that the conductor $N_L$ satisfies a Northcott property, thanks to the work of Deligne (see \cref{ex:conductor_l-adic}).

\subsection{\texorpdfstring{$L$}{L}-functions of pure motives}
\label{sec:Northcott_pure_motives}

Let $F$ be a number field. Then, we denote by $\mathcal{M}_F$ the category of pure motives for absolute Hodge cycles, which has been introduced by Deligne in \cite{Deligne_1979} (see also \cite[\S~6]{Deligne_Milne_1982}). 

Given a prime $\ell \in \mathbb{N}$, we have an $\ell$-adic realisation functor:
\begin{equation} \label{eq:ell_realisation}
    R_\ell \colon \mathcal{M}_F \to \operatorname{Rep}_{\mathbb{Q}_\ell}(G_F)
\end{equation}
which associates to every pure motive an $\ell$-adic representation of the absolute Galois group $G_F := \operatorname{Gal}(\overline{F}/F)$. 
Therefore, for every prime ideal $\mathfrak{p} \subseteq \mathcal{O}_F$, every prime $\ell \in \mathbb{N}$ and every motive $X \in \mathcal{M}_F$, one can define a polynomial
\begin{equation}
    P_\mathfrak{p}(R_\ell(X),t) := \det(1-\operatorname{Frob}_\mathfrak{p} t \mid D_\mathfrak{p}(R_\ell(X))) \in \mathbb{Q}_\ell[t]
\end{equation}
where $\operatorname{Frob}_\mathfrak{p} \subseteq G_F/I_\mathfrak{p}$ denotes the conjugacy class of geometric Frobenius elements at $\mathfrak{p}$, which is well defined only up to the action of the inertia subgroup $I_\mathfrak{p}$. 
To this end, when $\ell \nmid \lvert \mathcal{O}_F/\mathfrak{p} \rvert$, one lets these Frobenius elements act on the invariant sub-spaces $D_\mathfrak{p}(R_\ell(X)) := R_\ell(X)^{I_\mathfrak{p}} \subseteq R_\ell(X)$.
If on the other hand $\ell$ divides $\lvert \mathcal{O}_F/\mathfrak{p} \rvert$, the definition of $D_\mathfrak{p}(X)$ is more complicated, and requires $p$-adic Hodge theory, as explained in \cite[\S~3]{fo92}.
One conjectures that the polynomial $P_\mathfrak{p}(R_\ell(X),t)$ has always rational coefficients.
If this is true, then we can form the following formal Euler product
\begin{equation} \label{eq:motivic_L}
    L(R_\ell(X),s) := \prod_{\mathfrak{p}} \frac{1}{P_\mathfrak{p}(R_\ell(X),\lvert \mathcal{O}_F/\mathfrak{p}\rvert^{-s})}
\end{equation}
where we can plug in a complex variable $s \in \mathbb{C}$.

Now, each motive $X \in \mathcal{M}_F$ has a weight grading $X = \bigoplus_{w \in \mathbb{Z}} X_w$, as explained in \cite[Theorem~6.7]{Deligne_Milne_1982}.
Given $w \in \mathbb{Z}$, one says that $X$ is \textit{pure} of weight $w$ if $X = X_w$. Then, if $X \in \mathcal{M}_F$ is pure of weight $w$, one conjectures that for every prime $\ell \in \mathbb{N}$, the formal Euler-product $L(R_\ell(X),s)$ defined in \eqref{eq:motivic_L} gives rise to an element of $\mathcal{L}^{(w)}(F)$, with prescribed functional equation (see \cref{rmk:motivic_L} for details).
This motivates the following definition.

\begin{definition} \label{def:Mw}
    Given a number field $F$ and an integer $w \in \mathbb{Z}$, we denote by $\mathcal{M}^{(w)}(F)$ the set of isomorphism classes of motives $X \in \mathcal{M}_F$ such that:
    \begin{itemize}
        \item $X$ is pure of weight $w$;
        \item there exists a prime $\ell \in \mathbb{N}$ such that:
        \begin{itemize}
            \item for every prime ideal $\mathfrak{p} \subseteq \mathcal{O}_F$, we have that $P_{\mathfrak{p}}(R_\ell(X),t) \in \mathbb{Q}[t]$;
            \item the Euler product $L(R_\ell(X),s)$ defined in \eqref{eq:motivic_L} can be meromorphically continued to a function $L(R_\ell(X),s) \in \mathcal{L}^{(w)}(F)$ whose degree and conductor are given by: 
            \begin{equation} \label{eq:conductor_motivic_L}
                \begin{aligned} 
                d_{L(R_\ell(X),s)} &= \dim(R_\ell(X)) \\
                N_{L(R_\ell(X),s)} &= \lvert \Delta_F \rvert^{\dim(R_\ell(X))} \cdot \lvert \operatorname{N}_{F/\mathbb{Q}}(\mathfrak{f}_{R_\ell(X)}) \rvert
            \end{aligned}
            \end{equation}
            where $\Delta_F$ is the discriminant of the number field $F$ (see \cref{ex:discriminants}), and $\mathfrak{f}_{R_\ell(X)}$ is the conductor of the Galois representation $R_\ell(X)$ (see \cref{ex:conductor_l-adic}).
        \end{itemize}
    \end{itemize}
\end{definition}

\begin{remark}
    One also hopes that the polynomials $P_\mathfrak{p}(R_\ell(X),t)$ are independent of $\ell$.
    Therefore, one often writes $P_\mathfrak{p}(X,t)$ and $L(X,s)$ instead of $P_\mathfrak{p}(R_\ell(X),t)$ and $L(R_\ell(X),s)$.
    However, we will not need this in our paper, as the results we prove hold for every prime $\ell \in \mathbb{N}$.
\end{remark}
\begin{remark}
    In fact, 
    one can also give an intrinsic definition of the $L$-function $L(X,s)$ of a motive $X \in \mathcal{M}_F$, which does not depend on the choice of an auxiliary prime $\ell$, by taking the action of $\operatorname{Frob}_\mathfrak{p}$ on $D_\mathfrak{p}(R_p(X))$, where $p \in \mathbb{N}$ is the rational prime underlying $\mathfrak{p}$. This is the point of view adopted by Deninger in \cite{de94}.
\end{remark}

We are now ready to prove \cref{thm:Northcott_left_critical_strip}, which shows that special values of $L$-functions of pure motives of a fixed dimension have the Northcott property, if they are taken to the left of the critical strip.

\begin{proof}[Proof of {\cref{thm:Northcott_left_critical_strip}}]
    Since the $L$-functions associated to the $\ell$-adic realisation of any element of $\mathcal{M}^{(w)}(F)$ belongs to $\mathcal{L}^{(w)}(F)$ by assumption, \cref{prop:Selberg_bound} shows that there exists a constant $B_3$, depending only on $B_1, B_2, w$ and $n$, such that for every $X \in \mathcal{M}_{\mathbf{B}}^{(w)}(F,n)$ we have that $N_{L(R_\ell(X),s)} \leq B_3$. Moreover, we have the equality:
    \[
        N_{L(R_\ell(X),s)} = \lvert \Delta_F \rvert^{\dim(R_\ell(X))} \cdot \mathcal{C}_0(R_\ell(X))
    \]
    by definition of $\mathcal{M}^{(w)}(F)$, where $\mathcal{C}_0(R_\ell(X)) := \lvert \operatorname{N}_{F/\mathbb{Q}}(\mathfrak{f}_{R_\ell(X)}) \rvert$. Therefore, we can use the fact that the pair: 
    \[
        \{\dim,\mathcal{C}_0 \colon \mathcal{G}_\ell^{(w)}(F) \to \mathbb{R} \}    
    \]
    has the Northcott property (as explained in \cref{ex:conductor_l-adic}) to conclude that the set $R_\ell(S_{B_1,B_2})$ is finite. This allows us to conclude, because the realisation functor \eqref{eq:ell_realisation} is fully faithful (as proved by Jannsen in \cite[Theorem~4.4]{ja90}).
\end{proof}

\begin{remark} \label{rmk:motivic_L}
    We recall that the set $\mathcal{M}^{(w)}(F)$ is supposed to coincide with the set of isomorphism classes of motives $X \in \mathcal{M}_F$ which are pure of weight $w$. 
    In particular, for every rational prime $\ell \in \mathbb{N}$, the $L$-function $L(R_\ell(X),s)$ associated to a pure motive $X \in \mathcal{M}_F$ of weight $w$ is conjectured to belong to the set $\mathcal{L}^{(w)}(F)$ defined in \cref{def:L_class}. 
    
    Indeed, as we mentioned above, $L$-functions of motivic origin are usually defined by the Euler product \eqref{eq:motivic_L}, which is then conjectured to converge for $\Re(s) > w/2$. 
    The functional equation \eqref{eq:functional_equation} has also historically been conjectured to be true (as explained in Fontaine's \cite[\S~12.4]{fo92} and in Deninger's \cite[\S~3]{de94} works), with an explicit formula for all the invariants appearing in \eqref{eq:functional_equation}, as explained for instance by Serre in \cite[\S~4]{Serre_1969}. 
    In particular, the equalities portrayed in \eqref{eq:conductor_motivic_L} should always be true.

    To conclude, let us observe that the bound $\lvert \alpha_{j,\mathfrak{p}} \rvert \leq \operatorname{N}_{F/\mathbb{Q}}(\mathfrak{p})^{w/2}$, which appears in the definition of $\mathcal{L}^{(w)}(F)$, is explicitly mentioned in \cite[Footnote~32]{fo92}, and is a consequence of Deligne's weight monodromy conjecture.
    More precisely, let $V$ be a smooth and projective variety defined over a number field $F$, and let $X = H^w(V) \in \mathcal{M}_F$ for some $w \in \mathbb{N}$. 
    If $\mathfrak{p} \subseteq \mathcal{O}_F$ is a prime ideal at which $V$ has good reduction, and $\ell \in \mathbb{N}$ is any prime such that $\ell \nmid \lvert \mathcal{O}_F/\mathfrak{p} \rvert$, then one can combine the smooth and proper base change theorem in étale cohomology, together with the Weil conjectures (proven by Deligne), to show that the roots of $P_\mathfrak{p}(R_\ell(X),t)$ have absolute value equal to $\lvert \mathcal{O}_F/\mathfrak{p} \rvert^{w/2}$.
    On the other hand, one conjectures that if $\mathfrak{p}$ is a prime of bad reduction for $V$ then the absolute values of the roots of $P_\mathfrak{p}(R_\ell(X),t)$ are of the form $\lvert \mathcal{O}_F/\mathfrak{p} \rvert^{i/2}$, for some $i \in \{0,\dots,w\}$.

    More generally, one conjectures that if $K$ is a discretely valued field of residual characteristic $p$, and $\ell \in \mathbb{N} \setminus \{p\}$ is a prime number, then for every variety $V$ defined over $K$, the weights appearing in the $\ell$-adic cohomology groups $H^w_{\text{ét}}(V_{\overline{K}};\mathbb{Q}_\ell)$ should lie in the set $\{0,\dots,2w\}$, and in the smaller set $\{0,\dots,w\}$ if $V$ is proper (as noted already by Serre and Tate in \cite[Appendix, Problem~2]{Serre_Tate_1968}). 
    Similar properties are known to hold for the $\ell$-adic cohomology of varieties defined over finite fields (as proved by Deligne in \cite[Théorème~3.3.1]{Deligne_1980}), and for the singular cohomology of complex algebraic varieties (as shown by Deligne in \cite[Théorème~8.2.4]{Deligne_1974}).
    Moreover, the aforementioned conjecture follows from Deligne's weight monodromy conjecture, outlined in Illusie's survey \cite[\S~3.8]{Illusie_1994}. This conjecture is known to hold when $V$ is an abelian variety, thanks to the work of Grothendieck \cite{Grothendieck_1972}. This implies the validity of the weight monodromy conjecture for the weight $w = 1$, and any variety $V$. Moreover, Rapoport and Zink \cite{Rapoport_Zink_1982} verified this conjecture when $w = 2$, and the more recent work of Scholze \cite{Scholze_2012} has shown that the aforementioned conjecture holds true if $V$ can be realised as a set-theoretic complete intersection inside a toric variety.
\end{remark}

\begin{remark} \label{rmk:mixed_Northcott}
    It is not easy to dispose of the purity assumption appearing in \cref{thm:Northcott_left_critical_strip}.
    More precisely, let $F$ be a number field, and $\mathcal{MM}_F$ denote the category of mixed motives defined by Jannsen in \cite{ja90}, which contains Deligne's category $\mathcal{M}_F$ as a full subcategory (as explained in \cite[Theorem~4.4]{ja90}).
    Now, each object $X \in \mathcal{MM}_F$ admits a weight filtration $\mathcal{W}$ with finitely many graded pieces. 
    Hence, for every $X \in \mathcal{MM}_F$, there exist $w_{\min}(X), w_{\max}(X) \in \mathbb{Z}$ such that $\operatorname{gr}_w^\mathcal{W}(X) = 0$ if either $w < w_{\min(X)}$ or $w > w_{\max}(X)$.
    Moreover, to every $X \in \mathcal{MM}_F$ one can associate an $L$-function $L(X,s)$, a completed $L$-function $\widehat{L}(X,s)$ and a dual $X^\vee \in \mathcal{MM}_F$, such that $\widehat{L}(X,s)$ and $\widehat{L}(X^\vee,1-s)$ should be related by a functional equation, as explained in \cite{de94}.
    
    Now, for every $w \in \mathbb{Z}$, we denote by $\mathcal{MM}^{(w)}(F)$ the set of isomorphism classes of those $X \in \mathcal{MM}_F$ such that $w_{\min}(X) > w$, and with the property that $\widehat{L}(X,s)$ is a well defined meromorphic function which satisfies the expected functional equation.
    Moreover, for every $\mathbf{B} = (B_1,B_2,B_3) \in \mathbb{R}_{\geq 0}^3$ and every $w, n \in \mathbb{Z}$ with $2 n < w$, we define the set:
    \[
        \mathcal{MM}^{(w)}_{\mathbf{B}}(F,n) := \{ X \in \mathcal{MM}^{(w)}(F) \colon \dim(X) \leq B_1, \ \lvert L^\ast(X,n) \rvert \leq B_2, \ w_{\max}(X) \leq B_3 \}
    \]
    which is analogous to the set $\mathcal{M}^{(w)}_\mathbf{B}(F,n)$ appearing in \eqref{eq:finiteMM}. 
    Then, adapting the proofs of \cref{prop:Selberg_bound} and \cref{thm:Northcott_left_critical_strip}, one can see that the set $\{X^\text{ss} \colon X \in \mathcal{MM}^{(w)}_{\mathbf{B}}(F,n) \}$ is finite, where $X^\text{ss} := \bigoplus_{j \in \mathbb{Z}} \operatorname{gr}_j^\mathcal{W}(X)$ is the semi-simplification of $X$ with respect to the weight filtration.

    However, this is not enough to conclude that the set $\mathcal{MM}^{(w)}_{\mathbf{B}}(F,n)$ is finite. In fact, we expect this set to be always infinite, at least as soon as $B_3 - w \geq 1$. Indeed, as Deninger already remarks in \cite{de94}, the $L$-function $L(X,s)$ defined in \cite{de94} is not substantially different from the $L$-function $L(X^\text{ss},s)$. This is readily seen for the case of Kummer motives (described by Scholl in \cite[\S~2]{Scholl_1994}), which should provide the simplest class of mixed motives for which the special values at the left of the critical strip of the corresponding $L$-functions do not have the Northcott property. 
    We aim at coming back to this in future work.
\end{remark}

\begin{remark}
    Let $F$ be a number field, and $h \colon \mathcal{MM}_F \to \mathbb{R}$ denote Kato's height of mixed motives, which we introduced in \cref{ex:Heights_for_motives}.
    Then, Kato shows in \cite[Proposition~2.1.17]{ka18} that a Northcott property for this height would have momentous consequences, such as the finite generation of the ``motivic cohomology'' groups given by extensions in Jannsen's category $\mathcal{MM}_F$.
    More precisely, to prove such a finite generation one would need to show that $h$ satisfies a Northcott property when restricted to motives $X \in \mathcal{MM}_F$ with fixed semi-simplification. 
    As we outlined in the introduction, it seems promising to do so by relating Kato's height to special values of $L$-functions. Nevertheless, the previous \cref{rmk:mixed_Northcott} suggests that Deninger's notion of an $L$-function of mixed motives may not be well suited for this purpose, since it does not differentiate enough between a mixed motive and its semi-simplification. On the other hand, the special values of the multivariate $L$-function associated to a mixed motive $X$ by the pioneering work of Brown \cite{Brown_2019} may be more closely related to Kato's height. This will be the subject of future investigations.
\end{remark}

We conclude this section by applying \cref{thm:Northcott_left_critical_strip} to $L$-functions of elliptic curves and abelian varieties.

\begin{corollary} \label{cor:abelian_varieties_outside}
    Let $F$ be a number field and $g \in \mathbb{Z}_{\geq 1}$. We denote by $\mathcal{A}_g(F)$ the set of isomorphism classes of $g$-dimensional abelian varieties $A$ defined over $F$. Moreover, we denote by $\mathcal{A}'_g(F) \subseteq \mathcal{A}_g(F)$ the subset of those abelian varieties $A$ whose associated $L$-function $L(A,s)$ belongs to $\mathcal{L}^{(1)}(F)$. 
    Then, for every $n \in \mathbb{N}$ and $B \in \mathbb{R}_{\geq 0}$, the set:
    \begin{equation} \label{eq:abelian_set}
        \{ A \in \mathcal{A}'_g(F) \colon \lvert L^\ast(A,-n) \rvert \leq B \}
    \end{equation}
    is finite.
\end{corollary}
\begin{proof}
    Thanks to \cite[Theorem~6.25]{Deligne_Milne_1982} we have an embedding $\mathcal{A}_g'(F) \hookrightarrow \mathcal{M}^{(1)}(F)$, which sends $A$ to $H^1(A)$. Moreover, $\dim(H^1(A)) = 2 \dim(A)$, which yields an embedding of the set \eqref{eq:abelian_set} into the set $\mathcal{M}^{(1)}_{B,2g}(F)$. Therefore, the set \eqref{eq:abelian_set} is finite, thanks to \cref{thm:Northcott_left_critical_strip}.
\end{proof}

Now, it has been classically conjectured that $\mathcal{A}_g'(F) = \mathcal{A}_g(F)$. This is known when $g = 1$ and $F = \mathbb{Q}$, thanks to the celebrated modularity theorem, proved by Breuil, Conrad, Diamond and Taylor in \cite{Breuil_Conrad_Diamond_Taylor_2001}, extending ideas of Taylor and Wiles. More generally, we know that $\mathcal{A}_g'(F) = \mathcal{A}_g(F)$ when $g=1$ and $F$ is a totally real number field, thanks to the potential modularity result proven in Winterberger's appendix to \cite{Nekovar_2009}.
Finally, the recent potential modularity results of Boxer, Calegari, Gee and Pilloni \cite{Boxer_Calegari_Gee_Pilloni_2021} show that $\mathcal{A}_g'(F) = \mathcal{A}_g(F)$ when $g=2$ and $F$ is totally real, or $g = 1$ and $F$ is a quadratic extension of a totally real number field.
Therefore, we get the following completely unconditional corollary of \cref{thm:Northcott_left_critical_strip}.

\begin{corollary} \label{cor:modular_abelian_varieties_left}
    Fix a natural number $n \in \mathbb{N}$ and a real number $B \in \mathbb{R}_{\geq 0}$. Then, for every totally real number field $F$ and every extension $K \supseteq F$ such that $[K \colon F] \leq 2$, the sets:
    \[
        \begin{aligned}
            \{E \in \mathcal{A}_1(K) &\colon \lvert L^\ast(E,-n) \rvert \leq B \} \\
            \{A \in \mathcal{A}_2(F) &\colon \lvert L^\ast(A,-n) \rvert \leq B \}
        \end{aligned}
    \]
    are finite.
\end{corollary}

\subsection{Number fields}
\label{sec:Northcott_NF_left}

The aim of this section is to extend the results of \cref{thm:Northcott_left_critical_strip} to the simplest set of motives of fixed weight but unbounded dimension, which is given by the Weil restrictions $\mathbf{1}_{F/\mathbb{Q}} := \operatorname{N}_{F/\mathbb{Q}}(\mathbf{1}_F) \in \mathcal{M}^{(0)}(\mathbb{Q})$ of the trivial motives $\mathbf{1}_F := H^0(\operatorname{Spec}(F)) \in \mathcal{M}^{(0)}(F)$, where $F$ varies in the class of number fields. 
For these motives, one has that $d_F := \dim(\mathbf{1}_{F/\mathbb{Q}}) = [F \colon \mathbb{Q}]$ and $L(\mathbf{1}_{F/\mathbb{Q}},s) = \zeta_F(s)$, where $\zeta_F$ denotes Dedekind's $\zeta$-function.
Therefore, if we let $F$ vary, we see that the dimension of our motives is no longer bounded. Nevertheless, we can use the exponential growth of the discriminant $\lvert \Delta_F \rvert$ with respect to the degree $d_F$ in order to show the following effective result, which is part of \cref{thm:Northcott_class_number_formula}.

\begin{proposition} \label{prop:Dedekind_left}
    Let $S(\mathbb{Q})$ denote the set of isomorphism classes of number fields. Moreover, for every $n \in \mathbb{Z}$, we define the family of sets:
    \begin{equation} \label{eq:sets_of_number_fields}
        S_{B}(\mathbb{Q},n) := \{ [F] \in S(\mathbb{Q}) \colon \lvert \zeta^\ast_F(n) \rvert \leq B \}
    \end{equation}
    depending on $B \in \mathbb{R}_{\geq 0}$. 
    Then, we have that: 
    \begin{itemize}
        \item for every $n \in \mathbb{Z}_{\geq 2}$ and every $B \in \mathbb{R}$ such that $B \geq \zeta(n)^2$, where $\zeta(s) = \zeta_\mathbb{Q}(s)$ denotes Riemann's zeta function, the set $S_B(\mathbb{Q},n)$ is infinite;
        \item for every $n \in \mathbb{Z}_{\leq -1}$ and $B \in \mathbb{R}_{\geq 0}$, the set $S_B(\mathbb{Q},n)$ is finite.
    \end{itemize}
    Moreover, there exists an absolute, effectively computable constant $c_0 \in \mathbb{R}_{> 0}$ such that for every $B \in \mathbb{R}_{> 1}$ the following bound:
    \begin{equation} \label{eq:effective_left}
        \lvert S_{B}(\mathbb{Q},n) \rvert \leq \exp\left( \frac{c_0}{1 - n} \log\left( B \right) \left( \log\log\left( B \right) \right)^3 \right)
    \end{equation}
    holds true.
\end{proposition}
\begin{proof}
    First of all, we recall that for every $s \in \mathbb{C}$ such that $\Re(s) > 1$, the function $\zeta_F(s)$ admits various expansions as a Dirichlet series and an Euler product. More precisely, we have that
\begin{equation} \label{eq:Euler_product_Dedekind}
    \zeta_F(s) := \sum_{0 \neq I \subseteq \mathcal{O}_F} \frac{1}{\operatorname{N}_{F/\mathbb{Q}}(I)^s} = \sum_{n=1}^{+\infty} \frac{a_n(F)}{n^s} = \prod_{\mathfrak{p} \subseteq \mathcal{O}_F} \frac{1}{1 - \operatorname{N}_{F/\mathbb{Q}}(\mathfrak{p})^{-s}}
\end{equation}
where $a_n(F)$ denotes the number of ideals $I \subseteq \mathcal{O}_F$ such that $\operatorname{N}_{F/\mathbb{Q}}(I) := \lvert \mathcal{O}_F/I \rvert$.
In particular, we see from \eqref{eq:Euler_product_Dedekind} that: 
\begin{equation} \label{eq:trivial_bound_zeta}
    1 \leq \zeta_F(s) = \prod_{\mathfrak{p} \subseteq \mathcal{O}_F} \frac{1}{1 - \operatorname{N}_{F/\mathbb{Q}}(\mathfrak{p})^{-s}} \leq \prod_p \left( \frac{1}{1-p^{-s}} \right)^{a_p(F)} \leq \zeta(s)^{d_F}
\end{equation}
for every $s \in \mathbb{R}_{> 1}$. Therefore, we see that for every $n \in \mathbb{Z}_{\geq 2}$, the set $S_B(\mathbb{Q},n)$ is infinite as soon as $B \geq \zeta(n)^2$.

    Now, let us prove the second part of the statement. We will proceed by mimicking the proof of \cref{prop:Selberg_bound}. 
    First of all, let $r_1(F)$ denote the number of embeddings $F \hookrightarrow \mathbb{R}$, and $r_2(F)$ denote the number of complex conjugate pairs of embeddings $F \hookrightarrow \mathbb{C}$ whose image is not contained in $\mathbb{R}$. Then, the completed $\zeta$-function:
\[
    \widehat{\zeta}_F(s) := \widehat{L}(\mathbf{1}_{F/\mathbb{Q}},s) = \Gamma_\mathbb{R}(s)^{r_1(F)} \Gamma_\mathbb{C}(s)^{r_2(F)} \zeta_F(s)
\] 
satisfies the following functional equation:
\begin{equation} \label{eq:functional_equation_Dedekind}
    \widehat{\zeta}_F(s) = \lvert \Delta_F \rvert^{\frac{1}{2} - s} \widehat{\zeta}_F(1-s)
\end{equation}
for every $s \in \mathbb{C}$ (see \cite[Chapter~VII, Corollary~5.10]{ne99}).

    Therefore, we see that for every $[F] \in S_B(\mathbb{Q},n)$ the following upper bound:
    \begin{equation} \label{eq:Dedekind_left_1}
        \lvert \Delta_F \rvert^{\frac{1}{2} - n} \leq \left\lvert \frac{\Gamma_\mathbb{R}^\ast(n)}{\Gamma_\mathbb{R}^\ast(1-n)} \right\rvert^{r_1(F)} \cdot \left\lvert \frac{\Gamma_\mathbb{C}^\ast(n)}{\Gamma_\mathbb{C}^\ast(1-n)} \right\rvert^{r_2(F)} \cdot \frac{B}{\zeta_F(1-n)}
    \end{equation}
    holds true. We would like to use such an upper bound to conclude that $\lvert \Delta_F \rvert$ is bounded, and thus that $S_B(\mathbb{Q},n)$ is finite, thanks to Hermite's theorem. This will also yield explicit upper bounds for the cardinality of $S_B(\mathbb{Q},n)$, thanks to Couveignes's result which we recalled in \cref{ex:discriminants}.

     To show that $\lvert \Delta_F \rvert$ is indeed bounded, for every $[F] \in S_B(\mathbb{Q},n)$, we apply first of all \eqref{eq:trivial_bound_zeta} to \eqref{eq:Dedekind_left_1}, which yields the following upper bound:
    \[
        \lvert \Delta_F \rvert^{\frac{1}{2} - n} \leq \left\lvert \frac{\Gamma_\mathbb{R}^\ast(n)}{\Gamma_\mathbb{R}^\ast(1-n)} \right\rvert^{r_1(F)} \cdot \left\lvert \frac{\Gamma_\mathbb{C}^\ast(n)}{\Gamma_\mathbb{C}^\ast(1-n)} \right\rvert^{r_2(F)} \cdot B
    \]
    for every $[F] \in S_B(\mathbb{Q},n)$.
    Now, let: 
    \[
        \gamma_\mathbb{R}(n) := \left\lvert \frac{\Gamma_\mathbb{R}^\ast(n)}{\Gamma_\mathbb{R}^\ast(1-n)} \right\rvert
    \]
    and observe that:
    \[
        \gamma_\mathbb{R}(n) = \begin{cases}
            \frac{(2 \pi)^{2 m}}{(2 m)!}, \ &\text{if} \ n = - 2m \\
            \pi \frac{(2 \pi)^{2 m}}{(2 m)!}, \ &\text{if} \ n = - 2m - 1
        \end{cases}
    \]
    for every $n \in \mathbb{Z}_{\leq -1}$. Analogously, we have that:
    \[
        \gamma_\mathbb{C}(n) := \left\lvert \frac{\Gamma_\mathbb{C}^\ast(n)}{\Gamma_\mathbb{C}^\ast(1-n)} \right\rvert = (2 \pi) \cdot \left( \frac{(2 \pi)^{-n}}{(-n)!} \right)^2
    \]
    for every $n \in \mathbb{Z}_{\leq -1}$.
    Therefore, one sees easily that:
    \[
        \gamma_\mathbb{R}(n) \leq \alpha_1 := \gamma_\mathbb{R}(-7) = \pi \cdot \frac{(2 \pi)^6}{6!}
    \]
    for every $n \in \mathbb{Z}_{\leq -1}$. Analogously, we have that:
    \[
        \gamma_\mathbb{C}(n) \leq \alpha_2 := \gamma_\mathbb{C}(-6) = 2 \cdot \pi \cdot \left(\frac{(2 \pi)^6}{6!}\right)^2
    \]
    for every $n \in \mathbb{Z}_{\leq -1}$. Therefore, the following upper bound holds true:
    \begin{equation} \label{eq:discriminant_upper_bound}
        \lvert \Delta_F \rvert^{\frac{1}{2} - n} \leq \gamma_\mathbb{R}(n)^{r_1(F)} \cdot \gamma_\mathbb{C}(n)^{r_2(F)} \cdot B \leq \alpha_1^{r_1(F)} \cdot \alpha_2^{r_2(F)} \cdot B
    \end{equation}
    for every $[F] \in S_B(\mathbb{Q},n)$.

    Now, recall that there exists an effectively computable constant $d_0 \in \mathbb{N}$ such that:
    \begin{equation} \label{eq:poitou}
        \lvert \Delta_F \rvert \geq (60.8)^{r_1(F)} (497.2)^{r_2(F)}
    \end{equation}
    for every number field $F$ such that $d_F \geq d_0$. Indeed, this is one of the most optimised unconditional versions of Stark's lower bounds for discriminants, which is due to Poitou (see \cite[Equation~(K)]{Poitou_1977}).
    Therefore, if we set:
    \[
        T_B(\mathbb{Q},n) := \{ [F] \in S_B(\mathbb{Q},n) \colon d_F \leq d_0 \} 
    \]
    then we know that $T_B(\mathbb{Q},n)$ is finite, thanks to \cref{thm:Northcott_left_critical_strip}. In fact, we see directly from \eqref{eq:discriminant_upper_bound} that the following upper bound:
    \begin{equation} \label{eq:disc_bound_small_degree}
        \lvert \Delta_F \rvert \leq \left( \alpha_1^{d_0} \cdot B \right)^{\frac{2}{1-2n}}
    \end{equation}
    holds true for every $[F] \in T_B(\mathbb{Q},n)$. This implies that $T_B(\mathbb{Q},n)$ is finite, thanks to Hermite's theorem. To conclude, define:
    \begin{equation*} 
        g_n := \max\left( \frac{\log(\gamma_\mathbb{R}(n))}{\log(60.8)}, \frac{\log(\gamma_\mathbb{C}(n))}{\log(497.2)},0 \right)
    \end{equation*}
    and let us observe that: 
    \begin{equation} \label{eq:bound_gn}
        1-2(n+g_n) \geq (1 - 2(g_{-1}-1))(-n) = (1.223\dots) \cdot (-n)
    \end{equation}
    for every $n \in \mathbb{Z}_{\geq -1}$. Hence, for every $[F] \in S_B(\mathbb{Q},n) \setminus T_B(\mathbb{Q},n)$ we have that:
    \begin{equation} \label{eq:disc_bound_large_degree}
        \lvert \Delta_F \rvert \leq B^{\frac{2}{1-2(n+g_n)}} \leq B^{\frac{1.64}{-n}}
    \end{equation}
    as follows by combining the bounds \eqref{eq:discriminant_upper_bound}, \eqref{eq:poitou} and \eqref{eq:bound_gn}.
    Therefore, we see from \eqref{eq:disc_bound_small_degree} and \eqref{eq:disc_bound_large_degree} that there exists an absolute constant $a_1 \in \mathbb{R}_{> 0}$ such that the following upper bound:
    \[
        \lvert \Delta_F \rvert \leq a_1 \cdot B^{\frac{1.64}{-n}}
    \]
    holds true for every $[F] \in S_{B}(\mathbb{Q},n)$. To conclude, we can simply apply \eqref{eq:couveignes} to see that \eqref{eq:effective_left} holds true.
\end{proof}

\begin{remark}
    In principle, the proof of \cref{prop:Dedekind_left} could be adapted to show that for every number field $F$ and every $w,n \in \mathbb{Z}$ such that $2 n < w$, the function: 
    \[
        \begin{aligned}
            \mathcal{M}^{(w)}(F) &\to \mathbb{R} \\
        [X] &\mapsto \lvert L^\ast(X,n) \rvert
        \end{aligned}
    \]
    has the Northcott property on its own, without having to bound the dimension $\dim(X)$.
    In order to do so, we would need to show the existence of two constants $C_w \in \mathbb{R}_{>1}$ and $D_w \in \mathbb{R}_{> 0}$, depending at most on $w$ and on the number field $F$, such that 
    \[
        \operatorname{N}_{F/\mathbb{Q}}(\mathfrak{f}_X) \geq C_w^{\dim(X)}
    \]
    for every $X \in \mathcal{M}^{(w)}(F)$ with $\dim(X) \geq D_w$.
    Such a lower bound, which would be an analogue of \eqref{eq:poitou} for more general motives, is known to hold (conditionally) for Artin $L$-functions, thanks to work of Odlyzko \cite{Odlyzko_1977}, and would be expected for abelian varieties. 
    Moreover, a similar lower bound, which is valid for every motive $X \in \mathcal{M}^{(w)}(F)$, has been proved by Mestre in \cite[\S~III]{Mestre_1986}.
    However, let us note that Mestre's lower bound does not guarantee that $C_w > 1$ when $w$ is big. More precisely, this can be achieved only when one can show that the absolute value of the sum of the roots of the polynomials $P_\mathfrak{p}(R_\ell(X),t)$ is much smaller than the trivial bound. This can be achieved for instance when $F = \mathbb{Q}$ and $X = H^1(A)$ for some abelian variety $A$ whose associated $L$-functions satisfies all the properties conjectured in \cref{def:L_class}. In this case, Mestre obtains in \cite[Page 229, Proposition]{Mestre_1986} that $\lvert \mathfrak{f}_A \rvert > 10^{\dim(A)}$.  
\end{remark}

\section{Special values inside the critical strip}\label{zero}

The aim of this section is to study the Northcott properties satisfied by the special values of $L$-functions taken inside the critical strip. We will concentrate on two particular cases:
\begin{itemize}
    \item Dedekind $\zeta$-functions $\zeta_F(s)$ associated to number fields $F$. In this case, we will show that the special value $\zeta_F^\ast(0)$ satisfies a Northcott property, whereas $\zeta_F^\ast(1)$ does not.
    In other words, the different behaviour of the values at the left and the right of the critical strip expressed in \cref{prop:Dedekind_left} holds true also at the boundary of the critical strip;
    \item $L$-functions $L(A,s)$ associated to abelian varieties $A$.
    In this case, we have only one special value $L^\ast(A,1)$ at the centre of the critical strip, which is the subject of the notorious conjecture of Birch and Swinnerton-Dyer.
    Even assuming this conjecture, we will see that it is surprisingly difficult to prove or even predict whether or not the special value $L^\ast(A,1)$ satisfies a Northcott property.
\end{itemize}

\subsection{Number fields}
\label{sec:boundary_NF}

In this section, we look at the special values $\zeta_F^\ast(0)$ and $\zeta_F^\ast(1)$ associated to a number field $F$. 
First of all, we show in \cref{prop:Northcott_class_number_formula} that the special value at $s = 0$ has the Northcott property, whereas the special value at $s = 1$ does not. Moreover, \cref{prop:Northcott_class_number_formula_quantitative} gives the quantitative result, with the proof of the explicit upper bound for the cardinality of the finite sets appearing in \cref{prop:Northcott_class_number_formula}.
Together with \cref{prop:Dedekind_left}, and its effective version given by \eqref{eq:effective_left}, these results complete the proof of \cref{thm:Northcott_class_number_formula}.

\begin{proposition}
\label{prop:Northcott_class_number_formula}
    For every $B \in \mathbb{R}_{> 0}$ the set $S_B(\mathbb{Q},0)$ defined in \eqref{eq:sets_of_number_fields} is finite. 
    On the other hand, there exists $B_1 \in \mathbb{R}_{> 0}$ such that for every $B \geq B_1$ the set $S_B(\mathbb{Q},1)$ is infinite.
\end{proposition}
\begin{proof}
    Fix a number field $F$.
    First of all, let us recall that the analytic class number formula (see \cite[Chapter~VII, Corollary~5.11]{ne99}) gives: 
    \begin{equation} \label{eq:CNF_1}
        \zeta_F^\ast(1) = \frac{2^{r_1(F)} \cdot (2 \pi)^{r_2(F)}}{w_F} \cdot \frac{h_F \cdot R_F}{\sqrt{\lvert \Delta_F \rvert}}
    \end{equation}
    where $h_F, R_F, w_F$ are respectively the class number, the regulator and the number of roots of unity of $F$.
    Combining this with the functional equation \eqref{eq:functional_equation_Dedekind}, we have the following expression:
    \begin{equation}\label{eq:CNF}
    \zeta_F^{\ast}(0) = - \frac{h_F}{w_F} \, R_F
    \end{equation}
    for the special value of $\zeta_F(s)$ at $s = 0$.
    In particular, \eqref{eq:CNF} shows that any upper bound for $\lvert \zeta_F^\ast(0) \rvert$ entails an upper bound for the product $h_F \, R_F$ and the degree $d_F := [F \colon \mathbb{Q}]$.
    Indeed, it is known that $R_F$ grows exponentially in $d_F$ (as proved by Zimmert \cite[Satz~3]{zi81}), whereas $w_F \leq 4 d_F^2$, as follows from the easy lower bound $2 \varphi(n) \geq \sqrt{n}$ for Euler's totient function $\varphi$.
    This shows that $w_F$ grows at most logarithmically in $R_F$, and thus that $h_F R_F$ and $d_F$ are bounded whenever $\lvert \zeta_F^\ast(0) \rvert$ is bounded.
    Finally, we can simply apply the Brauer-Siegel theorem (see \cite[Theorem~2]{Brauer_1947}) to show that the discriminant $\Delta_F$ is bounded from above whenever $\lvert \zeta_F^\ast(0) \rvert$ is. Hence, Hermite's theorem shows that the set $S_B(\mathbb{Q},0)$ is finite.
    
    Now, let us prove that the function $[F] \mapsto \lvert \zeta_F^\ast(1) \rvert$ does not have the Northcott property, following a suggestion by Asbjørn Christian Nordentoft. 
    Let $\mathcal{D} \subseteq \mathbb{Z}_{< 0}$ denote the set of discriminants of imaginary quadratic fields.
    Then, there exist $c_1, c_2 \in \mathbb{R}_{> 0}$ such that:
    \[
        \frac{1}{X} \sum_{\substack{D \in \mathcal{D} \\ \lvert D \rvert \leq X}} h_{\mathbb{Q}(\sqrt{D})} = c_1 \sqrt{X} (1 + O(e^{- c_2 \sqrt{\log(X)}}))
    \]
    when $X \to +\infty$, as explained in \cite[Theorem~5.1]{Barban_1966}. 
    Moreover, \eqref{eq:CNF_1} implies that
    \[
    \zeta^\ast_{\mathbb{Q}(\sqrt{D})}(1) = \pi \cdot \frac{h_{\mathbb{Q}(\sqrt{D})}}{\sqrt{\lvert \Delta_{\mathbb{Q}(\sqrt{D})}} \rvert}
    \]
    for every $D \in \mathcal{D}$ such that $D < -4$.
    Thus, we see from \eqref{eq:CNF_1} that there exists an infinite subset $\mathcal{D}' \subseteq \mathcal{D}$ and a real number $B_1 \in \mathbb{R}_{\geq 0}$ such that $\lvert \zeta_{\mathbb{Q}(\sqrt{D})}^\ast(1) \rvert \leq B_1$ for every $D \in \mathcal{D}'$.
    Therefore, the sets $S_B(\mathbb{Q},1)$ are infinite as soon as $B \geq B_1$, as we wanted to show.
\end{proof}

It is an interesting problem to find explicit upper bounds for the cardinality of the sets $S_B(\mathbb{Q},0)$ introduced in \eqref{eq:sets_of_number_fields},
whose finiteness is proved in \cref{prop:Northcott_class_number_formula}.
We note that the proof that we gave cannot lead to such explicit upper bounds, due to the non-effectivity of the general version of the theorem of Brauer and Siegel.
Nevertheless, if one restricts to CM fields, the Brauer-Siegel theorem can be made effective, as shown by Stark in \cite{st74}. Combining this with some effective lower bounds for the regulator $R_F$, one obtains the following result.

\begin{proposition}
\label{prop:Northcott_class_number_formula_quantitative}
    There exist two absolute, effectively computable constants $c_1, c_2 \in \mathbb{R}_{> 0}$ such that the following upper bound:
    \begin{equation} \label{eq:effective_0}
        \lvert S_B(\mathbb{Q},0) \rvert \leq \exp\left( c_1 B^{c_2 \log\log(B)} (\log\log(B))^3 \right)
    \end{equation}
    holds true for every $B \in \mathbb{R}_{> 1}$, where $\lvert S_B(\mathbb{Q},0) \rvert$ denotes the cardinality of the finite set $S_B(\mathbb{Q},0)$ defined in \eqref{eq:sets_of_number_fields}.
\end{proposition}
\begin{proof}
    First of all, let us derive an explicit upper bound for the regulator $R_F$ and the degree $d_F := [F \colon \mathbb{Q}]$ of every number field $[F] \in S_B(\mathbb{Q},0)$.
    To do so, we use the chain of inequalities:
    \begin{equation} \label{eq:effective_1}
        R_F \overset{(a)}{\leq} h_F R_F \overset{(b)}{\leq} w_F B \overset{(c)}{\leq} 4 d_F^2 B \overset{(d)}{\leq} 4 B \left(\frac{\log(R_F/c_3)}{\log c_4}\right)^2
    \end{equation}
    where  $c_3 := (11.5)^{-39}$ and $c_4 := 1.15$ are two absolute constants, $h_F := \lvert \operatorname{Pic}(\mathcal{O}_F) \rvert$ denotes the class number of $\mathcal{O}_F$ and $w_F := \lvert (\mathcal{O}_F^\times)_\text{tors} \rvert$ denotes the number of roots of unity of $\mathcal{O}_F$. 
    Here the inequality $(a)$ follows from the fact that $h_F \in \mathbb{Z}_{\geq 1}$; the inequality $(b)$ follows from the class number formula \eqref{eq:CNF}, because $[F] \in S_B(\mathbb{Q},0)$; the inequality $(c)$ comes from the simple bound $w_F \leq 4 d_F^2$; finally, the inequality $(d)$ comes from the bound
    \begin{equation} \label{eq:effective_2}
        c_3 c_4^{d_F} \leq R_F    
    \end{equation}
    which holds for every number field $F$ (as shown by Zimmert in \cite[Satz~3]{zi81}).
    Thus, \eqref{eq:effective_1}  implies that
    \begin{equation}
    \label{eq:effective_22}
    \frac{R_F}{(\log R_F)^2} \leq c_5 B
    \end{equation} for some absolute, effectively computable constant $c_5>0$. 
    Since there exists an absolute, effectively computable constant $c_6 > 0$ such that $x \geq c_6 \log(x)^3$ for every $x \in \mathbb{R}_{>0}$, we see that \eqref{eq:effective_22} implies that $\log R_F\leq \frac{c_6}{c_5} B$. Therefore, using this inequality in \eqref{eq:effective_22} we get that
    \begin{align}
        \label{eq:effective_3} R_F &\leq c_7 B (\log B)^2,
    \end{align}
for some absolute, effectively computable constant $c_7 \in \mathbb{R}_{> 0}$.
Moreover, \eqref{eq:effective_2} implies that
        \begin{align}
        \label{eq:effective_4} d_F &\leq c_8 \log B
    \end{align}
    for some absolute, effectively computable constant $c_8 \in \mathbb{R}_{> 0}$.
    
    Now, we want to use \eqref{eq:effective_3} and \eqref{eq:effective_4} to find an explicit upper bound for the discriminant $\lvert \Delta_F \rvert$ of every number field $[F] \in S_B(\mathbb{Q},0)$. In order to do so, we need to consider separately the case of CM fields, \textit{i.e.} of those totally imaginary number fields which are quadratic extensions of a totally real number field.
    More precisely, let $S^\text{CM}$ be the set of isomorphism classes of CM fields.
    Then, a combination of \cite[Theorem~C]{fr89} and \cite[Proposition~3.7]{pa14} shows that for every number field $[F] \in S_B(\mathbb{Q},0) \setminus S^\text{CM}$ one has the explicit upper bound:
    \begin{equation} \label{eq:effective_5}
        \lvert \Delta_F \rvert \leq d_F^{d_F} \, \exp\left( \max\left( 1, \frac{d_F^{2 d_F} R_F}{c_{9}} \right) \right)
    \end{equation}
    where $c_{9} \in \mathbb{R}_{> 0}$ is some absolute, effectively computable constant.
    Moreover, for every CM field $[F] \in S_B(\mathbb{Q},0) \cap S^\text{CM}$, one obtains, for every $\varepsilon \in (0,\frac{1}{2}]$, the upper bound:
    \begin{equation} \label{eq:effective_6}
        \lvert \Delta_F \rvert \leq \left( \frac{h_F \, d_F !}{c_{11}(\varepsilon)^{d_F}} \left( \left( \frac{d_F}{2} \right)^{\frac{d_F}{2}} \exp\left( \max\left( 1, \frac{d_F^{d_F} R_F}{2^{d_F} c_{10}} \right) \right) \right)^{\frac{1}{2} - \frac{1}{d_F}} \right)^{\frac{d_F}{d_F - 3 - \varepsilon d_F}}
    \end{equation}
    by using \eqref{eq:effective_5} for the maximal, totally real sub-field of $F$, and combining this estimate with \cite[Theorem~2]{st74} and \cite[Proposition~3.7]{pa14}.
    Here, $c_{11} \colon (0,\frac{1}{2}] \to \mathbb{R}$ is some effectively computable function, appearing in \cite{st74}.
    Finally, combining \eqref{eq:effective_5} and \eqref{eq:effective_6} (where we may fix $\varepsilon = 1/4$, for example) with the bounds \eqref{eq:effective_3} and \eqref{eq:effective_4}, we obtain the bound:
    \begin{equation} \label{eq:effective_7}
        \lvert \Delta_F \rvert \leq c_{12} \exp{({B^{c_{13}\log\log B}})}
    \end{equation}
    for every $[F] \in S_B$, and some effectively computable constants $c_{12}, c_{13} \in \mathbb{R}_{> 0}$. 
    
    Now, to conclude, we obtain \eqref{eq:effective_0} by summing \eqref{eq:couveignes} over $1 \leq d \leq c_9 \log(B)$ and replacing $X$ by $c_{12} \exp( B^{c_{13} \log\log(B)})$, which can be done thanks to \eqref{eq:effective_4} and \eqref{eq:effective_7}.
\end{proof}

\begin{remark} \label{rmk:Lalin}
    The results summarized in \cref{thm:Northcott_class_number_formula}, and proved in \cref{prop:Dedekind_left,prop:Northcott_class_number_formula,prop:Northcott_class_number_formula_quantitative}, have been generalized in two recent works by Généreux, Lalín and Li \cite{Genereux_Lalin_Li_2022}, and by Généreux and Lalín \cite{Genereux_Lalin_2022}. In particular, the former deals with the sets:
    \[
        S_B(\mathbb{F}_q(T),s) := \{ [F] \in S(\mathbb{F}_q(T)) \colon \lvert \zeta_F^\ast(s) \rvert \leq B \}
    \]
    where $q \in \mathbb{N}$ is a power of a prime number, and $S(\mathbb{F}_q(T))$ denotes the set of isomorphism classes of finite extensions of the global function field $\mathbb{F}_q(T)$.
    Moreover, if $q \equiv 1 (4)$ and $s \in \mathbb{C}$ is any complex number with $\Re(s) > \frac{1}{2}$, then \cite[Theorem~1.3]{Genereux_Lalin_Li_2022} shows that $S_B(\mathbb{F}_q(T),s)$ is infinite, for $B$ sufficiently large. On the other hand, the same theorem shows that for every prime power $q \in \mathbb{N}$, every $B \in \mathbb{R}_{\geq 0}$ and every $s \in \mathbb{C}$ with $\Re(s) < \frac{1}{2} - \frac{\log(2)}{\log(q)}$, the set $S_B(\mathbb{F}_q(T),s)$ is finite, and its cardinality can be explicitly bounded (see \cite[Theorem~3.6]{Genereux_Lalin_Li_2022}). On the other hand, \cite{Genereux_Lalin_2022} deals with the sets $S_B(\mathbb{Q},s)$ defined in \eqref{eq:sets_of_number_fields}, where one replaces the integer $n \in \mathbb{Z}$ by any complex number $s \in \mathbb{C}$.
    In particular, they show that the presence of poles of the $\Gamma$-function at negative integers implies that there are regions on the left of the critical strip where the Northcott property does not hold (see \cite[\S~4.3]{Genereux_Lalin_2022}). 
    Moreover, \cite[Theorem~3.2]{Genereux_Lalin_2022} and \cite[Theorem~6.6]{Genereux_Lalin_2022} prove that the Northcott property does not hold for $\Re(s) \geq \frac{1}{2}$, which includes in particular the right half of the critical strip. These phenomenons had also been noticed by the authors during private communications with Jerson Caro. 
\end{remark}

\subsection{Abelian varieties}
\label{sec:abelian_varieties}

In this section, we investigate the possible Northcott property of the special value at the integer $s = 1$ of the $L$-functions $L(A,s) := L(H^1(A),s)$ associated to abelian varieties $A$ defined over a number field $F$. The main outcome of the discussion is that even assuming the Birch and Swinnerton-Dyer conjecture, it is not possible, so far, to prove a Northcott property in this case. According to the heuristics of Watkins about elliptic curves $E$ over $\mathbb{Q}$ (see \cite{Watkins_2008}), one is in fact led to the conclusion that the Northcott property for $L^*(E,1)$ could be unlikely to hold in general.

We note that it is not clear whether we can expect a similar Northcott property as in the case of the special values $\zeta_F^\ast(0)$ which we considered in the previous section.
First of all, if we want to follow the strategy that we used in the previous section, we should relate the special value $L^\ast(A,1)$ to some regulator determinant.
This relation was given by the class number formula \eqref{eq:CNF} in the case of the special value $\zeta^\ast_F(0)$ studied in the previous section, and was thus unconditional.
On the other hand, the celebrated conjecture by Birch and Swinnerton-Dyer (see (1.5) page 418 and (B) page 419 of \cite{ta66}) predicts that $L^\ast(A,1)$ is related to a regulator determinant $R_{A/F}$ by the equality
\begin{equation} \label{eq:bsd_formula_new}
    L^{\ast}(A,1) \stackrel{?}{=} c_F \left( \frac{\prod_{v \in M^0_F} c_v(A)}{\lvert A(F)_{\text{tors}} \rvert \, \lvert A^{\vee}(F)_{\text{tors}}\rvert} \right)  \frac{\lvert \Sha(A/F) \rvert  R_{A/F}}{\Omega_{A}^{-1}},
\end{equation}
where $c_v(A)$ stands for the Tamagawa number at the place $v$, we denote $A^{\vee}$ for the dual abelian variety, $\Sha(A/F)$ is the Tate-Shafarevitch group, $\Omega_{A}$ is the archimedean period, and $c_F$ is a quantity depending on the discriminant and the degree of the number field $F$.

Now, the first step in the proof of \cref{thm:Northcott_class_number_formula} was observing that the quantity $w_F = \lvert (\mathcal{O}_F^\times)_\text{tors} \rvert$ appearing in the class number formula \eqref{eq:CNF} is clearly bounded from above by a polynomial in the degree $[F \colon \mathbb{Q}]$ of the number field $F$.
An analogous statement for abelian varieties is the content of the following conjecture.
\begin{conjecture}[Torsion conjecture] \label{conj:torsion}
    For every $d\in \mathbb{N}_{\geq 1}$ and every $g \in \mathbb{N}_{\geq 1}$ there exists a natural number $c(g,d) \in \mathbb{N}$ such that for all number fields $F$ of degree $d=[F:\mathbb{Q}]$, for all $g$-dimensional abelian varieties $A$ defined over $F$, we have $\lvert A(F)_{\text{tors}} \rvert \leq c(g,d)$.
\end{conjecture}

We recall that, in the case of elliptic curves, \cref{conj:torsion} is proved to be true, thanks to work of Merel (see \cite{me96}).
In higher dimension the conjecture is still open, but the prime number theorem shows easily that 
\[
\lvert A(F)_{\text{tors}} \rvert \, \lvert A^\vee(F)_{\text{tors}} \rvert \ll (\log\lvert \operatorname{N}_{F/\mathbb{Q}}(\mathfrak{f}_A) \rvert)^{4 \operatorname{dim}(A)}
\] 
as explained in \cite[Lemma~3.6]{hi07}.

Now for a fixed field of definition $F$, observing that the Tamagawa numbers $c_v(A)$ are positive integers, we see that any upper bound for the quantity $\lvert L^\ast(A,1) \rvert$ entails an upper bound for the quantity
\begin{equation} \label{eq:Northcott_BSD_ratio}
    \frac{\lvert \Sha(A/F) \rvert  R_{A/F}}{\Omega_{A}^{-1}}
\end{equation}
if one assumes the validity of the formula \eqref{eq:bsd_formula_new}, and of \cref{conj:torsion}.
Since our goal is to study the Northcott property for the quantity $\lvert L^\ast(A,1) \rvert$, it would be useful to compare the quantity \eqref{eq:Northcott_BSD_ratio} to other quantities for which a Northcott property is already known to hold.
The best candidates for this are the stable Faltings height $h_\text{st}(A)$ and the norm of the conductor ideal $\mathfrak{f}_A$ of the abelian variety $A$.

This is exactly the same strategy which we employed in the proof of \cref{thm:Northcott_class_number_formula}, where the quantity $h_F \, R_F$ was compared to the quantity $\lvert \Delta_F \rvert$, which satisfies the Northcott property thanks to Hermite's theorem.
However, there is one fundamental difference between the proof of \cref{thm:Northcott_class_number_formula} and the current discussion: both the numerator and the denominator of the ratio \eqref{eq:Northcott_BSD_ratio} are comparable to something satisfying a Northcott property, at least conjecturally. In that respect, that case of $\vert L^*(A,1)\vert$ is closer to the case of $\vert\zeta_F^*(1)\vert$ described at the end of the previous section.

Let us be more precise.
First of all, if one denotes by $H(A)$ the exponential of the stable Faltings height of $A$, one has that
\[
    H(A) \ll \Omega_A^{-1} \ll H(A) (\log(H(A)))^{\operatorname{dim}(A)/2}
\]
as shown in \cite[Lemma~3.7]{hi07}.
Recent works of Hindry \cite{hi07}, Hindry-Pacheco \cite{Hindry_Pacheco_2016} and Griffon \cite{Griffon_2018_LMS,Griffon_2018_JNT,Griffon_2018_JTNB,Griffon_2019} on the analogue of the Brauer-Siegel estimate for abelian varieties show that the numerator of \eqref{eq:Northcott_BSD_ratio} is also expected to be comparable (in some cases) to $H(A)$. 
Hence it is necessary to gain further evidence in order to be able to decide if a Northcott property for the special value $L^\ast(A,1)$ associated to abelian varieties holds in some cases.
In particular, both the numerator and denominator of the ratio \eqref{eq:Northcott_BSD_ratio} appear to be comparable (in some cases) to $H(A)$.
This makes it extremely difficult to prove, or even expect, a Northcott property for the special value $L^\ast(E,1)$.
In fact, the heuristics proposed by Watkins in \cite{Watkins_2008} provide some evidence to expect that $L^\ast(E,1)$ does not satisfy a Northcott property. 
Indeed, Watkins's work predicts the existence of infinitely many elliptic curves $E$ defined over $\mathbb{Q}$ for which $\lvert \Sha(E/F) \rvert \, R_{E/F}$ is bounded (see in particular \cite[\S~4.5]{Watkins_2008}).
 
In the case of elliptic curves, one knows from \cite{ahp18} that the following inequality holds
\begin{equation} \label{eq:AHP_inequality}
        \frac{R_{E/F}}{\lvert E(F)_{\text{tors}} \rvert \, \lvert E^{\vee}(F)_{\text{tors}} \rvert} \gg h^{\frac{r_{E/F} - 4}{3}}  \left( \log(3 h) \right)^{\frac{2 r_{E/F} + 2}{3}}
\end{equation}
where $r_{E/F} := \operatorname{rk}(E(F))$ and $h := \max\{1,h(j(E))\}$. Here $h(j(E))$ is the absolute logarithmic Weil height of the $j$-invariant $j(E) \in F$ (see \cref{ex:Weil_height}), which is comparable with the stable Faltings height of $E$ (see for instance \cite[Lemma~3.2]{Pazuki_2018}).
The inequality \eqref{eq:AHP_inequality} shows that a part of the right hand side of \eqref{eq:bsd_formula_new} can indeed be related to some height, even if this relation is too weak to conclude (even assuming the validity of the Birch and Swinnerton-Dyer conjecture) that the special value $L^\ast(E,1)$ satisfies a Northcott property.

\section*{Acknowledgements}
We thank Fran\c{c}ois Brunault, Jerson Caro, Marco d'Addezio, Richard Griffon, Roberto Gualdi, Marc Hindry, Matilde Lalín, Asbjørn Christian Nordentoft and Martin Widmer for useful discussions. 
We thank moreover the anonymous referees for their helpful comments and suggestions.

\section*{Funding}
The first author is supported by ANR-17-CE40-0012 Flair and ANR-20-CE40-0003 Jinvariant. 
The second author performed this work within the framework of the LABEX MILYON (ANR-10-LABX-0070) of the Université de Lyon, within the program ``Investissements d'Avenir'' (ANR-11-IDEX-0007) operated by the French National Research Agency (ANR). He is also thankful to the Max Planck Institute for Mathematics in Bonn for its hospitality and financial support.
Both authors thank the IRN GANDA for its support.

\vspace{\baselineskip}
\noindent
\framebox[\textwidth]{
\begin{tabular*}{0.96\textwidth}{@{\extracolsep{\fill} }cp{0.84\textwidth}}
\raisebox{-0.7\height}{%
    \begin{tikzpicture}[y=0.80pt, x=0.8pt, yscale=-1, inner sep=0pt, outer sep=0pt, 
    scale=0.12]
    \definecolor{c003399}{RGB}{0,51,153}
    \definecolor{cffcc00}{RGB}{255,204,0}
    \begin{scope}[shift={(0,-872.36218)}]
      \path[shift={(0,872.36218)},fill=c003399,nonzero rule] (0.0000,0.0000) rectangle (270.0000,180.0000);
      \foreach \myshift in 
           {(0,812.36218), (0,932.36218), 
    		(60.0,872.36218), (-60.0,872.36218), 
    		(30.0,820.36218), (-30.0,820.36218),
    		(30.0,924.36218), (-30.0,924.36218),
    		(-52.0,842.36218), (52.0,842.36218), 
    		(52.0,902.36218), (-52.0,902.36218)}
        \path[shift=\myshift,fill=cffcc00,nonzero rule] (135.0000,80.0000) -- (137.2453,86.9096) -- (144.5106,86.9098) -- (138.6330,91.1804) -- (140.8778,98.0902) -- (135.0000,93.8200) -- (129.1222,98.0902) -- (131.3670,91.1804) -- (125.4894,86.9098) -- (132.7547,86.9096) -- cycle;
    \end{scope}
    \end{tikzpicture}%
}
&
Riccardo Pengo received funding from the European Research Council (ERC) under the European Union’s Horizon 2020 research and innovation programme (grant agreement number 945714).
\end{tabular*}
}

\printbibliography

\end{document}